\documentclass[11pt]{amsart}
\usepackage[T1]{fontenc}
\usepackage{lmodern}
\usepackage{xcolor}
\usepackage{amsmath,amssymb,amsthm,mathtools}
\usepackage[margin=1in]{geometry}
\usepackage{microtype}
\usepackage{enumitem}
\usepackage[hidelinks]{hyperref}
\hypersetup{
  pdftitle={Large entropy implies H\"older regularity on large sets},
  pdfauthor={Fabrizio Bianchi and Yan Mary He},
  pdfkeywords={holomorphic endomorphisms, entropy, pluripolar sets}
}

\newtheorem{theorem}{Theorem}[section]
\newtheorem{proposition}[theorem]{Proposition}
\newtheorem{lemma}[theorem]{Lemma}
\newtheorem{corollary}[theorem]{Corollary}
\theoremstyle{remark}
\newtheorem{remark}[theorem]{Remark}
\numberwithin{equation}{section}

\newcommand{\C}{\mathbb C}

\newcommand{\PP}{\mathbb P}

\newcommand{\R}{\mathbb R}
\newcommand{\Z}{\mathbb Z}
\newcommand{\supp}{\operatorname{supp}}
\newcommand{\dist}{\operatorname{dist}}
\newcommand{\Lip}{\operatorname{Lip}}
\newcommand{\Jac}{\operatorname{Jac}_{\R}}
\newcommand{\PSH}{\operatorname{PSH}}
\newcommand{\htop}{h_{\mathrm{top}}}
\newcommand{\ddc}{dd^c}
\newcommand{\one}{\mathbf 1}
\newcommand{\Ucal}{\mathcal U}

\title[Large entropy and H\"older regularity]{Large entropy implies H\"older regularity on large sets}
\author{Fabrizio Bianchi}
\address{Dipartimento di Matematica\\
Universit\`a di Pisa\\
Largo Bruno Pontecorvo 5, 56127 Pisa, Italy}
\email{fabrizio.bianchi@unipi.it}
\author{Yan Mary He}
\address{Department of Mathematics\\
University of Oklahoma\\
Norman, OK 73019, USA}
\email{he@ou.edu}
\date{\today}

\begin{document}

\begin{abstract}
Let $F$ be a holomorphic endomorphism of $\PP^k(\C)$ of algebraic
degree $d\ge2$, with $k\ge1$. We prove that every ergodic
$F$-invariant probability measure $\mu$ of entropy strictly
greater than
$(k-1)\log d$ admits compact sets of $\mu$-measure arbitrarily
close to one on which the restrictions of $\mu$ have H\"older
continuous super-potentials.
Consequently, $\mu$ gives zero mass
to every Borel pluripolar set.
\end{abstract}

\maketitle

\section{Introduction}
We write $\PP^k$ for the complex projective
space $\PP^k(\C)$.
Let $F\colon\PP^k\to\PP^k$ be a holomorphic endomorphism of
algebraic degree $d\ge2$, with $k\ge1$. Its topological entropy is
$k\log d$ \cite{Gro03}, and it admits a unique measure $\mu_F$ of maximal entropy;
see \cite{FLM83,Lyu83,Man83} when $k=1$, and
\cite{FS95,BriDuv01,DS10} when $k\ge2$.
In every dimension, $\mu_F$ is a Monge--Amp\`ere measure \cite{BT82} 
of a
H\"older continuous potential; see for instance \cite{DS10}.
By \cite[Theorem 1.1]{DN14}, this is equivalent to $\mu_F$ having
a H\"older continuous super-potential.
In particular, $\mu_F$ is \emph{moderate}
and gives zero mass to
pluripolar sets (and, as a consequence, to proper analytic subsets); see \cite{DNS10,DN14}.
The super-potential of a measure is a functional
on positive closed $(1,1)$-currents
and can be
identified
with the integration functional on their normalized
quasi-plurisubharmonic potentials; see \cite{DS09SP}.
The notion of super-potential
generalizes to higher dimensions the classical
notion of potential for measures on $\PP^1$.

In any dimension,
endomorphisms of $\PP^k$ 
admit large families of measures whose entropy is arbitrarily close to the maximal value $k\log d$, see for instance
\cite{Dup12,UZ13,BD23,BD24}.
The above discussion naturally raises the question of whether
sufficiently large entropy forces
analogous regularity for invariant measures other than $\mu_F$. For instance, as a consequence of the variational principle and Gromov entropy bound argument \cite{Gro03}, 
every ergodic $F$-invariant probability measure of entropy strictly
larger than $(k-1)\log d$ gives zero mass to every proper analytic subset; see
\cite[Section 1.7]{DS10}. 
However,
global H\"older regularity cannot be expected, even for mixing
measures with entropy arbitrarily close to the maximum.
Our main result gives a Lusin-type conclusion: every ergodic
invariant measure of entropy
strictly greater than $(k-1)\log d$ has
H\"older continuous super-potentials after restriction to compact
sets of measure arbitrarily close to one. As a consequence, every such measure also gives no mass to pluripolar subsets.

\subsection{H\"older regularity on arbitrarily large sets}
Let $\omega$ be the Fubini--Study form on $\PP^k$, normalized by
$\int_{\PP^k}\omega^k=1$. We write $h_\mu(F)$ for the metric entropy of $F$
with respect to an $F$-invariant probability measure $\mu$, and
$\htop(F)$ for the topological entropy of $F$.

A \emph{quasi-plurisubharmonic function} is locally the sum of a
plurisubharmonic function and a smooth function, and is not
identically $-\infty$. Set $\ddc:=\frac{i}{\pi}\partial\bar\partial$.
We denote by $\PSH(\PP^k,\omega)$ the quasi-plurisubharmonic
functions $u$ for which $\omega+\ddc u$ is a positive current;
these functions are often called \emph{$\omega$-plurisubharmonic};
see for instance
\cite{GZ05} and \cite[Appendix A]{DS10}.
For $v\in L^1(\omega^k)$, set
$\|v\|_1:=\int_{\PP^k}|v|\,\omega^k$.
For a Borel set $A\subset\PP^k$, let $\one_A$ denote its indicator
function and $\one_A\mu$ the restriction of $\mu$ to $A$.
Our main result is the following.

\begin{theorem}
\label{thm:submeasure}
Let $F\colon\PP^k\to\PP^k$ be a holomorphic endomorphism of
algebraic degree $d\ge2$, with $k\ge1$, and let $\mu$ be an
ergodic $F$-invariant probability measure with
$h_\mu(F)>(k-1)\log d$. There exists
$\beta\in(0,1)$ such that, for every $0<\theta<1$, there are a
compact set $E_\theta\subset \PP^k$ with $\mu(E_\theta)>1-\theta$ and
a constant $C_\theta<\infty$ for which
$\eta_\theta:=\one_{E_\theta}\mu$ integrates every
quasi-plurisubharmonic function and satisfies
\begin{equation}\label{eq:holder-restriction}
 \int_{\PP^k}|u-v|\,d\eta_\theta
 \le C_\theta\max\{\|u-v\|_1,\|u-v\|_1^\beta\}
 \qquad\text{for every }
 u,v\in\PSH(\PP^k,\omega).
\end{equation}
\end{theorem}

As mentioned above, the functional $u\mapsto\int_{\PP^k}u\,d\eta_\theta$, acting on
$\omega$-plurisubharmonic functions normalized by
$\int_{\PP^k}u\,\omega^k=0$, represents the normalized
super-potential
of $\eta_\theta$ under the identification
$u\leftrightarrow\omega+\ddc u$;
see \cite{DS09SP} and
\cite[Definition 3.4]{DN14}.
Estimate \eqref{eq:holder-restriction} makes this functional
H\"older continuous with respect to the $L^1(\omega^k)$ distance;
see \cite[Lemma 3.5]{DN14}.
An explicit admissible range for $\beta$ is given in
Remark \ref{rem:holder-exponent}.
When $k=1$, the entropy assumption reduces to $h_\mu(F)>0$, and
the conclusion follows 
directly
from Ma\~n\'e's dimension formula
\cite[Theorem A]{Man88} and classical potential theory;
see
Section \ref{subsec:dimension-one}.

\smallskip

The restrictions in Theorem \ref{thm:submeasure} imply vanishing
on pluripolar sets and the following entropy bounds.
Recall that a set $P\subset\PP^k$ is \emph{locally pluripolar} if, near each
point, $P$ is contained in the pole set of a plurisubharmonic function
that is not identically $-\infty$ on the connected neighbourhood. By the global Josefson theorem, this is equivalent to the existence of $u\in\PSH(\PP^k,\omega)$ such that $P\subset\{u=-\infty\}$; see \cite{Ale81}.

\begin{corollary}
\label{cor:entropy}
Let $F\colon\PP^k\to\PP^k$ be a holomorphic endomorphism of
algebraic degree $d\ge2$, with $k\ge1$. The following hold.
\begin{enumerate}
[label=\textup{(\roman*)},leftmargin=2em]
\item Every ergodic $F$-invariant probability measure $\mu$ with
$h_\mu(F)>(k-1)\log d$ gives zero mass to every Borel pluripolar
set $P\subset \PP^k$.
\item If an $F$-invariant probability measure $\nu$ satisfies
$\nu(P)=1$ for a Borel pluripolar set $P\subset \PP^k$, then
$h_\nu(F)\le(k-1)\log d$.
In particular, every non-empty compact pluripolar set $K$ with
$F(K)\subset K$ satisfies $\htop(F|_K)\le(k-1)\log d$.
\end{enumerate}
\end{corollary}

In dimension one, Corollary \ref{cor:entropy} is already a
consequence of classical results; see also the observation following
the questions in \cite[Section 2]{dTV15}.
Indeed,
Ma\~n\'e's dimension
formula 
gives positive pointwise dimension almost
everywhere for every ergodic invariant measure of
positive entropy,
whereas polar sets have Hausdorff dimension zero
\cite[Theorem 5.9.6]{AG01}.
A related vanishing statement holds for holomorphic automorphisms
of compact K\"ahler surfaces and, in the non-compact setting, for
compactly supported ergodic invariant measures of positive entropy
of polynomial automorphisms of H\'enon type and invertible
horizontal-like maps in complex dimension two \cite{Duj04HL}. 
These invertible maps are outside the large-topological-degree
setting of Corollary \ref{cor:entropy}. Although the conclusion follows from standard
ingredients
and is already known for the measures of
maximal entropy in these settings
\cite{BS91,BLS93,Duj04HL,Can14}, we have not found it explicitly stated in the
literature.
We therefore
include a short proof, combining
the Ledrappier--Young theory with a pluripolarity criterion of Dujardin; see Theorem \ref{thm:surface}.

\smallskip

The H\"older estimate \eqref{eq:holder-restriction} also gives
exponential integrability and a Monge--Amp\`ere representation
for each restriction $\eta_\theta$.
Set $\Ucal:=\{u\in\PSH(\PP^k,\omega):\sup_{\PP^k}u=0\}$.
Every quasi-plurisubharmonic function can be normalized into
$\Ucal$ by positive scaling and addition of a constant.
A finite positive measure $\eta$ is \emph{moderate} if
$\sup_{u\in\Ucal}\int_{\PP^k}e^{-\alpha u}\,d\eta<\infty$
for some $\alpha>0$; see \cite{DNS10}.

\begin{corollary}
\label{cor:moderate}
The restrictions $\eta_\theta$ in Theorem \ref{thm:submeasure}
have the following properties.
\begin{enumerate}
[label=\textup{(\roman*)},leftmargin=2em]
\item There is $\alpha>0$, independent of $\theta$, such that
$\sup_{u\in\Ucal}\int_{\PP^k}e^{-\alpha u}\,d\eta_\theta<\infty$
for every $0<\theta<1$. In particular, every $\eta_\theta$ is moderate.
\item For each $\theta$, there is a H\"older continuous function
$\varphi_\theta\in\PSH(\PP^k,\omega)$ such that
$\eta_\theta=\eta_\theta(\PP^k)(\omega+\ddc\varphi_\theta)^k$.
The H\"older exponent of $\varphi_\theta$ can be chosen independently
of $\theta$.
\end{enumerate}
\end{corollary}

The passage to the restrictions $\eta_\theta$ in
Theorem \ref{thm:submeasure} and Corollary \ref{cor:moderate}
is necessary in general: even a mixing invariant measure with
entropy arbitrarily close to the maximum need not integrate every
quasi-plurisubharmonic function; see
Proposition \ref{prop:nonintegrable}.

\subsection{Entropy bounds and extensions}
We next bound the \emph{Bowen entropy} of pluripolar subsets of uniformly
expanding sets.
For $Z\subset\PP^k$, let $h_{\mathrm B}(F,Z)$ denote the Bowen entropy,
defined using covers by Bowen balls of variable lengths; see
\cite{FH12} and
Section \ref{sec:bowen}.
We say that $F$ is \emph{uniformly expanding} on a compact set $K$
in the ambient tangent bundle if there are constants
$c>0$ and $\rho>1$ such
that $\|D(F^n)_xv\|\ge c\rho^n\|v\|$ for every $x\in K$,
$v\in T_x\PP^k$, and $n\ge1$, with norms taken in the
Fubini--Study metric.

\begin{theorem}
\label{thm:bowen}
Let $F\colon\PP^k\to\PP^k$ be a holomorphic endomorphism of
algebraic degree $d\ge2$, with $k\ge1$.
Let $K\subset \PP^k$ be a non-empty compact set with $F(K)\subset K$
on which $F$ is uniformly expanding in the ambient tangent bundle.
For every Borel pluripolar set $P\subset \PP^k$, we have
$h_{\mathrm B}(F,P\cap K)\le(k-1)\log d.$
\end{theorem} 

When $P\cap K$ is compact and forward invariant, the bound follows
from Corollary \ref{cor:entropy} and the variational principle
for topological entropy.
For a general Borel set $P\cap K$, however, an invariant probability
measure supported on $P\cap K$ need not exist, so the inverse branch
distortion estimates for invariant measures cannot be applied directly.
Our argument only requires contracting holomorphic inverse branches
with suitable quantitative bounds; we assume uniform expansion
for simplicity to ensure these bounds.
In dimension one, the uniform expansion assumption
is unnecessary. In fact, $h_{\mathrm B}(F,P)=0$
for every Borel polar set $P\subset\mathbb P^1$;
see Section \ref{subsec:dimension-one}.

\smallskip

Theorems \ref{thm:submeasure} and \ref{thm:bowen} and Corollaries
\ref{cor:entropy} and \ref{cor:moderate} extend to 
polynomial-like maps
in any dimension \cite{DS03},
provided that $d_t>d_{k-1}^*(f)$,
and to
surjective
holomorphic self-maps of connected compact
K\"ahler manifolds of dimension $k\ge1$,
provided that
$\lambda_k(F)>\lambda_{k-1}(F)$; see Appendix \ref{sec:kahler}.
In dimension two, together with Gromov's entropy
bound, this gives a negative answer to a question
of de Th\'elin and Vigny
\cite[Section 2]{dTV15}.
A related question of de Th\'elin and Vigny concerning the
dimension of positive closed currents was studied, in a
directional form, by Dupont and Rogue \cite{DR20}.

\subsection{Proof strategy}
The key idea behind the proofs of Theorems \ref{thm:submeasure} and 
\ref{thm:bowen}
is a 
regularization argument, which was inspired by
 the work of
F. Rodriguez Hertz, J. Rodriguez Hertz, and W. Wu \cite{RHW26}.
They apply heat regularization to the indicator of the singular set
and control the resulting loss of flow invariance through quantitative
commutator estimates.
We adapt this idea to measures of large entropy,
using inverse branches in place of the flow,
as in our earlier work on joining rigidity for rational maps
\cite{BH26Joining}.
Since integrability against quasi-plurisubharmonic
functions is not known a priori for these measures,
we construct regularizations and recover integrability
for suitable restrictions by controlling the
regularization errors.

Theorem \ref{prop:regularization} implements this idea through a
regularization criterion that does not require invariance of the measure.
We construct approximating measures by transporting radial averaging
kernels along contracting holomorphic inverse branches
\cite{BDM08,BD19}.
To compare two consecutive approximations, we express both averages
on a common coordinate domain.
In the spirit of the Dinh--Sibony $dd^c$-method \cite{DS03,DS10},
we decompose the pulled-back test function into a difference of local
potentials of positive closed $(1,1)$-currents and a pluriharmonic function.
The pluriharmonic term cancels, leaving a bound in terms of the masses
of these currents and the densities of the averaging measures.

The entropy assumption 
comes into play
when these local comparisons are summed, giving exponential decay
of the comparison errors,
after absorbing the small exponential losses from the shrinking domains.
Combining this decay with an exponential density bound yields
integrability and H\"older regularity
on sets of arbitrarily large measure,
proving Theorem \ref{thm:submeasure}.
Corollaries \ref{cor:entropy} and \ref{cor:moderate}
then follow from Theorem \ref{thm:submeasure}
by standard arguments and known results from
pluripotential theory.

For Theorem \ref{thm:bowen}, uniform expansion
provides the required inverse branches, while the
dynamical Frostman lemma \cite{FH12} gives
the Bowen ball mass bounds.

\subsection*{Organization of the paper}
Section \ref{sec:preliminaries} collects preliminary estimates,
and Section \ref{sec:proof} establishes the key
regularization criterion.
Section \ref{sec:applications} proves Theorems
\ref{thm:submeasure} and \ref{thm:bowen} and Corollaries
\ref{cor:entropy} and \ref{cor:moderate},
gives a direct proof of
Theorem \ref{thm:submeasure} in dimension one,
and discusses examples.
Section \ref{sec:surface} treats positive-entropy measures in the
surface automorphism, H\'enon, and horizontal-like settings.
Finally, Appendix \ref{sec:kahler} treats the K\"ahler and polynomial-like
extensions.

\subsection*{Acknowledgments}
The authors would like to thank Jana Rodriguez Hertz for her talk on \cite{RHW26} at the {\it Rigidity Theory of Diffeomorphisms and Flows} Workshop
at the Brin Center in October 2026, which 
was the inspiration for this work.

This project has received funding from the
project ANR TIGerS (ANR-24-CE40-3604)
and from the MIUR Excellence Department Project awarded to the
Department of Mathematics of the University of Pisa,
CUP I57G22000700001.
The first author is affiliated to the GNSAGA group of INdAM.
The second author is supported by National Science Foundation
DMS-2554272 and by the Simons Foundation Travel Support for
Mathematicians program.

\subsection*{AI Declaration}
The authors used ChatGPT 6 to assist with literature searches, to discuss and check mathematical arguments, and to suggest revisions to the exposition. 
More precisely, inspired by the talk by J. Rodriguez Hertz, 
the second author initially asked ChatGPT 6 whether the method used in \cite{RHW26} could be adapted to prove Corollary \ref{cor:entropy}, using inverse branches in place of geodesic flows in a similar way as in our earlier work \cite{BH26Joining}. 
ChatGPT 6 then proposed a possible strategy based on a similar regularization procedure
that required certain pluripotential-theoretic estimates. 
The authors subsequently developed these
estimates and improved the argument
to obtain the stronger Theorem \ref{thm:submeasure}.

\section{Preliminaries}
\label{sec:preliminaries}
This section collects potential-theoretic estimates and uniform
inverse branch and entropy estimates.

\subsection{Potential-theoretic estimates}
\label{subsec:potentials}
In this paper, distances, balls, and derivative norms on $\PP^k$ refer to the
Fubini--Study metric; in $\C^k$, balls
and norms are Euclidean. We write $\Lip(G)$ for the Lipschitz
constant of a map $G$ with respect to these metrics.
For a holomorphic map $G\colon U\to V$ between open subsets of
$\PP^k$, its real Jacobian is $\Jac G=(G^*\omega^k)/\omega^k$.
We set $d^c:=\frac{i}{2\pi}(\bar\partial-\partial)$.

For a positive closed $(1,1)$-current $S$ on $\PP^k$, write
$\|S\|=\int_{\PP^k}S\wedge\omega^{k-1}$
for its mass
and denote its cohomology
class by $[S]$. On a domain $V\subset\C^k$, write instead
$\|S\|_V=\int_VS\wedge\omega_{\mathrm{eucl}}^{k-1}$, where
$\omega_{\mathrm{eucl}}=\ddc|z|^2$, and denote Lebesgue measure by $dV$.

Recall that we denote
$\Ucal:=\{u\in\PSH(\PP^k,\omega):\sup_{\PP^k}u=0\}$.
For $u\in\Ucal$, the current $T_u:=\omega+\ddc u$ is positive
closed and has mass one. Compactness of normalized
$\omega$-plurisubharmonic functions in $L^1(\omega^k)$ gives
$\sup_{u\in\Ucal}\int_{\PP^k}(-u)\,\omega^k<\infty$;
see for instance \cite[Appendix A]{DS10} and
\cite[Proposition 2.7]{GZ05}.

The following is a direct consequence of \cite[Proposition A.16]{DS10}.

\begin{lemma}
\label{lem:comparison}
Fix $k\ge1$ and $r>0$, and write $B_t:=B(0,t)\subset\C^k$.
Suppose that $\sigma,\tau$ are probability measures supported in
$\overline{B_{2r}}$, absolutely continuous with respect to $dV$,
with densities bounded by some $M>0$, $dV$-almost everywhere.
Assume that
\begin{equation}\label{eq:pluriharmonic-cancellation}
 \int H\,d\sigma=\int H\,d\tau
 \qquad\text{for every real pluriharmonic function }H\text{ on }B_{3r}.
\end{equation}
Let $\varphi\in L^1_{\mathrm{loc}}(B_{4r})$ be real-valued and satisfy
$\ddc\varphi=S^+-S^-$, where $S^\pm$ are positive closed
$(1,1)$-currents of finite mass on $B_{4r}$.
Then there exists a constant $C>0$, depending only on $k$, such that
$$
 \left|\int\varphi\,d(\sigma-\tau)\right|
 \le CMr^2\bigl(\|S^+\|_{B_{4r}}+\|S^-\|_{B_{4r}}\bigr).
$$
\end{lemma}

\begin{proof}
By the local potential estimate in \cite[Proposition A.16]{DS10},
after rescaling to $B_{4r}$, there exist
plurisubharmonic functions $v^\pm$ on $B_{3r}$ with
$\ddc v^\pm=S^\pm$ and
$\|v^\pm\|_{L^1(B_{3r})}\le Cr^2\|S^\pm\|_{B_{4r}}$.
Since $\ddc(\varphi-v^++v^-)=0$, 
Weyl's lemma implies that
the function $\varphi-v^++v^-$
agrees $dV$-almost everywhere with a pluriharmonic function on $B_{3r}$.
Absolute continuity and \eqref{eq:pluriharmonic-cancellation}
imply that the integrals of $\varphi-v^++v^-$ against $\sigma$ and
$\tau$ agree. The density bounds applied to $v^+-v^-$ give the conclusion.
\end{proof}

Fix $r_*>0$, $C_*>1$, and
holomorphic charts $\Phi_x\colon B(0,r_*)\subset\C^k\to\PP^k$
for $x\in\PP^k$, with $\Phi_x(0)=x$ and $(x,z)\mapsto\Phi_x(z)$
Borel measurable. Uniformly in $x$ and $z,z'\in B(0,r_*)$, we have
\begin{equation}
\label{eq-charts}
 C_*^{-1}|z-z'|\le\dist(\Phi_x(z),\Phi_x(z'))\le C_*|z-z'|
 \qquad \text{and}
 \qquad C_*^{-1}dV\le\Phi_x^*\omega^k\le C_*dV.
\end{equation}
Such charts can be chosen from a finite holomorphic atlas.

We will
use 
radial averaging measures
to obtain exact
pluriharmonic cancellation. Fix a non-negative radial function
$\rho\in C_c^\infty(B(0,1))$ with $\int\rho\,dV=1$.
For $0<t<r_*$, set
\begin{equation}
\label{eq-kxt}
\kappa_{x,t}:=(\Phi_x)_*\bigl(t^{-2k}\rho(z/t)\,dV(z)\bigr).
\end{equation}
For every Borel set $A\subset\PP^k$, the map
$(x,t)\mapsto\kappa_{x,t}(A)$ is Borel measurable; thus
$(\kappa_{x,t})$ is a Borel probability kernel. Moreover, we have
\begin{equation}
\label{eq:kernel-properties}
 \supp\kappa_{x,t}\subset B(x,C_*t),
 \qquad
 0\le\frac{d\kappa_{x,t}}{d\omega^k}\le Ct^{-2k},
 \qquad \text{and} \quad
 \int H\,d\kappa_{x,t}=H(x),
\end{equation}
where the constant
$C>0$ is independent of $x$ and $t$, and the last identity
holds for every function $H$ which is pluriharmonic near
$\Phi_x(\overline{B(0,t)})$, by the Euclidean mean value property
for $H\circ\Phi_x$.

The next lemma 
is a standard telescoping regularity argument, stated
in a dual form for measures. It balances
approximation errors against density bounds
to obtain integrability and an absolute H\"older estimate; compare for instance
the regularity arguments in \cite[Lemma 1.19]{DS10}
and \cite[Lemma 2.11]{BD23}.

\begin{lemma}\label{lem:holder-limit}
Fix an integer $N\ge1$. Let $(\nu_n)_{n\ge N}$ be finite positive
measures on $\PP^k$ of equal mass, converging weakly to $\nu$.
Suppose that $a,b,C>0$ satisfy
\begin{equation}\label{eq:hp-telescopic}
 \nu_n\le Ce^{bn}\omega^k
 \quad \text{and}
 \quad
 \left|\int_{\PP^k}u\,d\nu_{n+1}-\int_{\PP^k}u\,d\nu_n\right|
 \le Ce^{-an}
 \qquad 
 \text{for every }
 n\ge N,\quad u\in\Ucal.
\end{equation}
Then $\nu$ integrates every quasi-plurisubharmonic function.
More precisely, setting
$\beta:=a/(a+b)$, there is $C'>0$, depending only on
$a,b,C,N$, such that
$$
 \int_{\PP^k}|u-v|\,d\nu
 \le C'\max\{\|u-v\|_1,\|u-v\|_1^\beta\}
 \qquad\text{for every }u,v\in\PSH(\PP^k,\omega).
$$
\end{lemma}

\begin{proof}
We first compare integrals of smooth test functions against $\nu_n$
and the limit measure $\nu$.
Since the measures have equal mass, subtracting $\sup_{\PP^k}u$
does not change any successive difference of integrals. Thus the
second
hypothesis 
in \eqref{eq:hp-telescopic}
applies to every smooth
$u\in\PSH(\PP^k,\omega)$.
For $m>n$, write
$\nu_m-\nu_n=\sum_{\ell=n}^{m-1}(\nu_{\ell+1}-\nu_\ell)$.
Summing the 
bounds in \eqref{eq:hp-telescopic}
and then letting $m\to\infty$ by weak
convergence gives, with $C_0:=C/(1-e^{-a})$,
\begin{equation}\label{eq:limit-comparison}
 \left|\int u\,d(\nu-\nu_n)\right|
 \le \sum_{\ell=n}^{\infty}
 \left|\int u\,d(\nu_{\ell+1}-\nu_\ell)\right|
 \le C_0e^{-an}.
\end{equation}

We next prove integrability and extend \eqref{eq:limit-comparison}
to every $u\in\PSH(\PP^k,\omega)$.
Choose smooth functions $u_j\in\PSH(\PP^k,\omega)$
decreasing to $u$; see, for instance,
\cite[Proposition A.41]{DS10}.
The density bound gives $u\in L^1(\nu_N)$, and
\eqref{eq:limit-comparison} applied to $u_j$ yields
$$\int u_j\,d\nu\ge\int u\,d\nu_N-C_0e^{-aN}.$$
The integrals $\int u_j\,d\nu$ therefore have a finite uniform lower
bound. Applying the monotone convergence theorem to the sequence
$u_1-u_j$ gives
$u\in L^1(\nu)$.
For each fixed $n\ge N$, the density bound also gives
$u\in L^1(\nu_n)$; passing to the limit in
\eqref{eq:limit-comparison} for $u_j$ proves that bound for $u$.
Scaling gives integrability of every quasi-plurisubharmonic function.

Next, 
we obtain a
H\"older bound for 
the integral of $u-v$
by balancing the
approximation error against the density bound.
For $u,v\in\PSH(\PP^k,\omega)$ and $t:=\|u-v\|_1 >0$,
\eqref{eq:limit-comparison} applied to $u$ and $v$, together with the first assumption in \eqref{eq:hp-telescopic}, 
gives
\begin{equation}\label{eq:pre-balance}
\left|\int(u-v)\,d\nu\right|\le C_1(e^{-an}+e^{bn}t).
\end{equation}
Choosing $n=\lceil\log(1/t)/(a+b)\rceil\ge N$ for sufficiently
small $t>0$, and $n=N$ for larger $t$, yields
\begin{equation}\label{eq:signed-bound}
\left|\int(u-v)\,d\nu\right|\le C_2\max\{t,t^\beta\}.
\end{equation}
The constants $C_1,C_2$ depend only on $a,b,C,N$.

We finally pass from the integral of $u-v$
to the integral of $|u-v|$ as in the statement.
Set $w:=\max(u,v)$, which belongs to $\PSH(\PP^k,\omega)$.
Since the functions $w-u$ and $w-v$ are non-negative, we have
$|u-v|=(w-u)+(w-v)$, and both $\|w-u\|_1$ and $\|w-v\|_1$
are at most $t$.
Applying
\eqref{eq:signed-bound} 
to the pairs $(w,u)$ and $(w,v)$ and
adding gives the stated estimate with $C'=2C_2$.
\end{proof}

\subsection{Inverse branches and entropy estimates}
\label{subsec:branches}
Let $\mu$ be an ergodic $F$-invariant probability measure on
$\PP^k$. To record inverse orbits, we use its \emph{natural extension}
$(\widehat X,\widehat F,\widehat\mu)$, where
$$
 \widehat X:=\{(x_j)_{j\in\Z}\in (\PP^k)^{\Z}:F(x_j)=x_{j+1}\},
 \qquad
 \widehat F((x_j)_{j\in\Z}):=(x_{j+1})_{j\in\Z}.
$$
For a history $\hat x=(x_j)_{j\in\Z}\in\widehat X$, write
$\pi_0(\hat x):=x_0$, so $F\circ\pi_0=\pi_0\circ\widehat F$.
The measure $\widehat\mu$ is the unique $\widehat F$-invariant
probability measure satisfying $(\pi_0)_*\widehat\mu=\mu$;
see \cite[Chapter 10]{CFS82}.

For $x\in \PP^k$, $\delta>0$, and $n\ge1$, let
$B_n(x,\delta):=\{y\in \PP^k:\dist(F^jy,F^jx)<\delta\text{ for }0\le j<n\}$
be the \emph{Bowen ball of length $n$ and radius $\delta$}.
A set is \emph{$(n,\delta)$-separated} if every two distinct points have
distance at least $\delta$ under some iterate $F^j$ with $0\le j<n$.

If $h_\mu(F)>(k-1)\log d$, then all Lyapunov exponents of
$\mu$ are strictly positive.
For $k=1$, Ruelle's entropy inequality gives a strictly positive Lyapunov
exponent whenever $h_\mu(F)>0$.
For $k\ge2$, positivity follows from \cite[Theorem C]{Dup12}, which
extends de Th\'elin's entropy--Lyapunov estimate \cite{dT08}
without assuming integrability of $\log\Jac F$.
In either case, positivity rules out the value $-\infty$, while
boundedness of $DF$ rules out $+\infty$, so all Lyapunov exponents
are finite.
Consequently, the multiplicative ergodic theorem gives a finite
limit for $n^{-1}\log\Jac(F^n)$ at $\mu$-almost every point.
Since $\log\Jac F$ is bounded above, the Birkhoff ergodic theorem applied to
its bounded lower truncations gives $\log\Jac F\in L^1(\mu)$:
otherwise,
$n^{-1}\log\Jac(F^n)$ would tend to $-\infty$
$\mu$-almost everywhere.

The next lemma combines the inverse branch estimates of
\cite{BDM08,BD19} with the Brin--Katok entropy formula
\cite{BK83} on one set of histories of large measure.

\begin{lemma}\label{lem:good-set}
Suppose that $h_\mu(F)>(k-1)\log d$, and fix
$(k-1)\log d<s<h_\mu(F)$ and $0<\theta<1$.
Let $\chi_{\min}>0$ denote the smallest Lyapunov exponent of $\mu$.
For every $0<\varepsilon<\chi_{\min}/2$, there exist a Borel
set $\widehat E\subset\widehat X$, constants
$r_0,\delta,\lambda>0$, $A\ge1$, and an integer $N\ge1$ such that
$\widehat\mu(\widehat E)>1-\theta$ and the following hold for every
$\hat x\in\widehat E$.
\begin{enumerate}[label=\textup{(\roman*)},leftmargin=2em]
\item For every $n\ge1$, there is an inverse branch
$g_{\hat x,n}\colon B(x_n,r_0e^{-\varepsilon n})\to \PP^k$
with $F^n\circ g_{\hat x,n}=\operatorname{id}$ and
$g_{\hat x,n}(x_n)=x_0$.
For all $y\in B(x_n,r_0e^{-\varepsilon n})$ and $0\le j\le n$, we have
$
 \dist(F^jg_{\hat x,n}(y),x_j)
 \le Ae^{\varepsilon n}e^{-\lambda(n-j)}\dist(y,x_n).
$
\item For every $n\ge N$, we have $\mu(B_n(x_0,\delta))\le e^{-sn}$.
\end{enumerate}
\end{lemma}

\begin{proof}
For \textup{(i)}, we use the tempered inverse branch estimates of
\cite[Theorem 1.4]{BDM08} and \cite{BD19}; see also
\cite[Section 2.2]{BH24}.
Fix $\lambda:=\chi_{\min}/2$.
The radius and distortion functions 
given by these estimates
are $\widehat\mu$-measurable, possibly with respect to the completed
$\sigma$-algebra. 
Since $\widehat X$ is a standard Borel space, they
admit Borel versions. Choose a Borel set $\widehat Y_0$ of full
$\widehat\mu$-measure on which these versions agree with the original
functions and all the inverse branch and tempered estimates hold,
and set
$\widehat Y:=\bigcap_{q\in\mathbb Z}
 \widehat F^q(\widehat Y_0)$.
Then $\widehat Y$ is an invariant Borel set of full
$\widehat\mu$-measure. On $\widehat Y$ we obtain Borel functions
$r\colon\widehat Y\to(0,1]$
and $L\colon\widehat Y\to[1,\infty)$ and, for
$\hat y=(y_j)_{j\in\Z}\in\widehat Y$, inverse branches
$b_{\hat y,m}$ of $F^m$ sending $y_0$ to $y_{-m}$, defined on
$B(y_0,r(\hat y))$ and satisfying
$\Lip(b_{\hat y,m})\le L(\hat y)e^{-\lambda m}$ for $m\ge1$.
The functions $r$ and $L$ can be chosen so that
$r(\widehat F^q\hat y)\ge e^{-\varepsilon q}r(\hat y)$ and
$L(\widehat F^q\hat y)\le e^{\varepsilon q}L(\hat y)$ for every integer $q\ge0$.

Choose $r_0>0$ and $A\ge1$ so that
$\widehat E_1:=\{\hat y\in\widehat Y:r(\hat y)\ge r_0,\ L(\hat y)\le A\}$
satisfies  $\widehat\mu(\widehat E_1)>1-\theta/2$.
For $\hat x\in\widehat E_1$, we have
$r(\widehat F^n\hat x)\ge e^{-\varepsilon n}r(\hat x)
\ge r_0e^{-\varepsilon n}$.
Thus we can define $g_{\hat x,n}$ as the restriction
of $b_{\widehat F^n\hat x,n}$ to $B(x_n,r_0e^{-\varepsilon n})$.
For $0\le j<n$, uniqueness of inverse germs gives
$F^j\circ g_{\hat x,n}=b_{\widehat F^n\hat x,n-j}$ on this ball.
The bound $L(\widehat F^n\hat x)\le Ae^{\varepsilon n}$ now gives
\textup{(i)}; the case $j=n$ follows from
$F^n\circ g_{\hat x,n}=\operatorname{id}$.

For \textup{(ii)}, the Brin--Katok entropy formula \cite{BK83}
and $s<h_\mu(F)$ give $\delta>0$, $N\ge1$, and a Borel set
$E_0\subset\PP^k$ with $\mu(E_0)>1-\theta/2$ such that
$\mu(B_n(x,\delta))\le e^{-sn}$ for every $x\in E_0$ and $n\ge N$.
Set $\widehat E:=\widehat E_1\cap\pi_0^{-1}(E_0)$.
Both conclusions hold on $\widehat E$, and
$\widehat\mu(\widehat E)>1-\theta$ since $(\pi_0)_*\widehat\mu=\mu$.
\end{proof}

\section{An inverse branch regularization criterion}
\label{sec:proof}
The following criterion deduces H\"older regularity
from the existence of contracting inverse branches
and Bowen ball mass bounds.
In the application to Theorem \ref{thm:submeasure},
$\Omega$ will be the good set of histories
$\widehat E$ provided by Lemma \ref{lem:good-set}.
The measure
$\eta$ is not
assumed to be invariant. When $\xi$ is fixed, we write $x_j$ for
$x_j(\xi)$.

\begin{theorem}
\label{prop:regularization}
Let $F\colon\PP^k\to\PP^k$ be a holomorphic endomorphism of
algebraic degree $d\ge2$, with $k\ge1$.
Let $m$ be a finite positive measure on a standard Borel space
$\Omega$, with $0<m(\Omega)\le1$, let $p\colon\Omega\to\PP^k$ be
Borel, and set $\eta:=p_*m$ and $x_j(\xi):=F^j(p(\xi))$ for
$\xi\in\Omega$ and $j\ge0$.
Fix $s>(k-1)\log d$ and $\varepsilon>0$ with
$(8k-4)\varepsilon<s-(k-1)\log d$.
Suppose that there are $r_0,\lambda,\delta,C_F>0$, $A\ge1$,
and an integer $N\ge1$ such that, for every $\xi\in\Omega$, the following hold.
\begin{enumerate}[label=\textup{(\roman*)},leftmargin=2em]
\item For every $n\ge1$,
there is a holomorphic map
$g_{\xi,n}\colon B(x_n,r_0e^{-\varepsilon n})\to\PP^k$
such that
$F^n\circ g_{\xi,n}=\operatorname{id}$,
 $g_{\xi,n}(x_n)=x_0$,
 and
$$
\dist(F^jg_{\xi,n}(y),x_j)
 \le Ae^{\varepsilon n}e^{-\lambda(n-j)}\dist(y,x_n),
 \qquad 0\le j\le n,\quad y\in B(x_n,r_0e^{-\varepsilon n}).
$$
\item For every $n\ge N$, we have $\eta(B_n(x_0,\delta))\le C_Fe^{-sn}$.
\end{enumerate}
Then $\eta$ integrates every quasi-plurisubharmonic function,
and there exist $\beta\in(0,1)$ and $C>0$ such that
\begin{equation}\label{eq:regularization-holder:revision}
 \int_{\PP^k}|u-v|\,d\eta
 \le C\max\{\|u-v\|_1,\|u-v\|_1^\beta\}
 \qquad\text{for all }u,v\in\PSH(\PP^k,\omega).
\end{equation}
The exponent $\beta$ can be chosen to depend only on
$s,\varepsilon,F$.
In particular, $\eta$ has a H\"older continuous super-potential.
\end{theorem}

Observe that no measurability of the family $(g_{\xi,n})_{\xi\in\Omega}$ is assumed.

\begin{proof}
Recall that we denote
$\Ucal=\{u\in\PSH(\PP^k,\omega):\sup_{\PP^k}u=0\}$;
see Section \ref{subsec:potentials}.
We will apply Lemma \ref{lem:holder-limit} by constructing measures
$\eta_n$ of mass $m(\Omega)$ that converge weakly to $\eta$ and
satisfy the two estimates in
\eqref{eq:hp-telescopic} for suitable constants $a,b,C>0$.

We use the charts $\Phi_x\colon B(0,r_*)\to\PP^k$ and their
metric and volume comparison constant $C_*$ fixed in
Section \ref{subsec:potentials}.
Throughout the proof, $C$ denotes a positive constant that may 
change from line to line,
but is independent of $n,j,\xi$, and $u\in\Ucal$.

\smallskip
\noindent{\bf Step 1. Constructing regularizations and grouping branches.}
The goal of this step is to construct absolutely continuous
approximations $\eta_n$ of $\eta$ that converge weakly to $\eta$
and satisfy the first bound in \eqref{eq:hp-telescopic}.
For each $n\ge N$,
we also construct a finite Borel partition of $\Omega$
into sets of $m$-mass $O(e^{-sn})$, on each of which
the corresponding restricted inverse branch is fixed.
For each $\xi$, we will replace the point mass at $x_0(\xi)$
by a radial averaging measure centred at $x_n(\xi)$ and transported back near $x_0(\xi)$
by $g_{\xi,n}$.
To verify that the transported radial measures can be integrated over $\Omega$,
we will
first restrict the branches $g_{\xi,n}$ to finitely many fixed domains.

\smallskip
\noindent\emph{Step 1.1. A Borel partition adapted to the inverse branches.}
Choose $a>0$ such that 
$$4a<r_*, \qquad 5C_*a<r_0/2,\qquad 
2C_*Aa<\delta$$
and set 
$R_n:=ae^{-2\varepsilon n}$.
For each $n$, choose a maximal $R_n/(2C_*)$-separated set of
centres $c_{n,j}$ in $\PP^k$, and set
$U_{n,j}^{(\ell)}:=\Phi_{c_{n,j}}(B(0,\ell R_n))$ for
$1\le\ell\le4$ and $U_{n,j}:=U_{n,j}^{(4)}$.
By maximality,
the balls $B(c_{n,j},R_n/(2C_*))$ cover $\PP^k$.
Since each such ball is contained in $U_{n,j}^{(1)}$, the sets
$U_{n,j}^{(1)}$ also cover $\PP^k$.
Moreover, $U_{n,j}\subset B(c_{n,j},4C_*R_n)$.
Hence, if
$x\in U_{n,j}$, then $c_{n,j}\in B(x,4C_*R_n)$. Since the centres
are $R_n/(2C_*)$-separated, the balls
$B(c_{n,j},{R_n}/{4C_*})$
are pairwise disjoint. For all indices $j$ such that
$x\in U_{n,j}$, these balls are contained in $B(x,CR_n)$, where
$C$ is independent of $n$ and $x$. Uniform volume comparison for
small balls on $\PP^k$ therefore shows that every point belongs
to at most a fixed number of the sets $U_{n,j}$, independently
of $n$.
Order the centres and set
$D_{n,j}:=U_{n,j}^{(1)}\setminus\bigcup_{i<j}U_{n,i}^{(1)}$.
Then $\{D_{n,j}\}_j$ is a finite Borel partition of $\PP^k$.

\smallskip
\noindent\emph{Step 1.2. Measurability of the regularizations.}
For $\xi$ with $x_n(\xi)\in D_{n,j}$,
the estimates for $\Phi_{c_{n,j}}$ in \eqref{eq-charts}
give
$\overline{U_{n,j}}
 \subset\overline{B(x_n(\xi),5C_*R_n)}
 \subset B(x_n(\xi),r_0e^{-\varepsilon n})$.
Thus every such $g_{\xi,n}$ is holomorphic on a neighbourhood
of $\overline{U_{n,j}}$ and restricts to an inverse branch of
$F^n$ on the same connected domain $U_{n,j}$.

Choose $t_0>0$ such that $t_0<r_*$ and
$2AC_*^2t_0<a$, and set
$$
 t_n:=t_0e^{-4\varepsilon n}
 \qquad
\text{and}
 \qquad
 K_{\xi,n}:=\Phi_{x_n}(B(0,t_n)).
$$
If $x_n\in D_{n,j}$, the chart comparisons give
$\dist(x_n,\partial U_{n,j}^{(2)})\ge R_n/C_*$ and
 $\overline{K_{\xi,n}}
 \subset\overline{B(x_n,C_*t_n)}$.
Since $C_*t_n<R_n/(2C_*)$, we have
$\overline{K_{\xi,n}}\subset U_{n,j}^{(2)}$.
In particular, $g_{\xi,n}$ is defined on a neighbourhood of
$\supp\kappa_{x_n,t_n}$
(see \eqref{eq-kxt}), and we set
$$
 Q_{\xi,n}:=(g_{\xi,n})_*\kappa_{x_n,t_n}.
$$
For each fixed $n$, we now partition $\Omega$
according to the index $j$ such that
$x_n(\xi)\in D_{n,j}$ and the restriction
of $g_{\xi,n}$ to $U_{n,j}$.
Let $\mathcal G_{n,j}$ be the collection of
restrictions $g_{\xi,n}|_{U_{n,j}}$ with
$x_n(\xi)\in D_{n,j}$ and, for $g\in\mathcal G_{n,j}$, set
$$
 \Omega_{n,j,g}:=
 \{\xi\in\Omega:
 x_n(\xi)\in D_{n,j},\ x_0(\xi)=g(x_n(\xi))\}.
$$
For fixed $n$ and $j$, choose a regular value
$y\in U_{n,j}$ of $F^n$. Evaluation at $y$ defines an injective
map
$$
\mathcal G_{n,j}\to (F^n)^{-1}(y),
\qquad g\mapsto g(y).
$$
Indeed, if two inverse branches agree at $y$, uniqueness of the
local inverse germ and connectedness of $U_{n,j}$ imply that they
agree identically. Since the fibre $(F^n)^{-1}(y)$ is finite,
$\mathcal G_{n,j}$ is finite.

If $\xi\in\Omega_{n,j,g}$, then
$g_{\xi,n}|_{U_{n,j}}$ and $g$ are inverse branches of $F^n$
on $U_{n,j}$ with the same value at $x_n(\xi)$. Hence, we have
$g_{\xi,n}|_{U_{n,j}}=g$.
Moreover, every $\Omega_{n,j,g}$ is Borel, since
$\xi\mapsto x_0(\xi)$ and $\xi\mapsto x_n(\xi)$ are Borel and
$g$ is continuous. 
Every $\xi\in\Omega$ belongs to one such set: the partition
$\{D_{n,j}\}_j$ determines $j$, and
$g_{\xi,n}|_{U_{n,j}}\in\mathcal G_{n,j}$ gives the corresponding
branch $g$. If $\xi$ belongs to both $\Omega_{n,j,g}$ and
$\Omega_{n,j,h}$, then $g(x_n)=h(x_n)$; uniqueness of the inverse
germ and connectedness of $U_{n,j}$ give $g=h$.
Thus, for every $n$,
the sets $\{\Omega_{n,j,g}\}_{j,g}$ form a finite Borel
partition of $\Omega$.

On $\Omega_{n,j,g}$, we have
$Q_{\xi,n}=g_*\kappa_{x_n(\xi),t_n}$.
In particular, to check measurability, observe that
for every Borel set $A\subset\PP^k$
we have
$
 Q_{\xi,n}(A)
 =\kappa_{x_n(\xi),t_n}(g^{-1}(A))$
 on $\Omega_{n,j,g}$.
The Borel dependence of the radial kernels on their centres and
the finiteness of the partition show that
$(Q_{\xi,n})_{\xi\in\Omega}$ is a Borel probability kernel.
We may therefore set
$$
 \eta_n:=\int_\Omega Q_{\xi,n}\,dm(\xi).
$$

\smallskip
\noindent\emph{Step 1.3. Convergence and density
bounds.}
We now prove weak convergence $\eta_n \to \eta$ and bound the densities of $\eta_n$.
Hypothesis \textup{(i)}, 
applied with $j=0$, and
\eqref{eq:kernel-properties} give, uniformly in $\xi$,
$\supp Q_{\xi,n}
 \subset B(x_0,Ce^{-(\lambda+3\varepsilon)n})$.
Take 
$\varphi\in C(\PP^k)$ and denote by $m_\varphi$
its modulus of continuity. Since
$\eta=p_*m$, we have
$$
 \left|\int_{\PP^k}\varphi\,d\eta_n
 -\int_{\PP^k}\varphi\,d\eta\right|
 \le
 \int_\Omega\int_{\PP^k}
 |\varphi(y)-\varphi(x_0(\xi))|
 \,dQ_{\xi,n}(y)\,dm(\xi)
 \le
 m(\Omega)\,
  m_\varphi\!\left(Ce^{-(\lambda+3\varepsilon)n}\right),
$$
which tends to zero as $n\to\infty$. Hence
$\eta_n\to\eta$ weakly. Moreover, all the measures $\eta_n$
have mass $m(\Omega)$.

Write
$\kappa_{x_n,t_n}=q_{\xi,n}\,\omega^k$, with
 $0\le q_{\xi,n}\le Ct_n^{-2k}$.
Since $F^n\circ g_{\xi,n}=\operatorname{id}$, changing variables
along the inverse branch gives
\begin{equation}
\label{eq:inverse-branch-density}
 \frac{dQ_{\xi,n}}{d\omega^k}(z)
 =
 \begin{cases}
 q_{\xi,n}(F^n(z))\Jac(F^n)(z),
 & z\in g_{\xi,n}(U_{n,j}),\\
 0, & \text{otherwise},
 \end{cases}
\end{equation}
where $j$ is determined by $x_n\in D_{n,j}$. Hence,
we have
$Q_{\xi,n}
 \le Ct_n^{-2k}\|\Jac(F^n)\|_\infty\,\omega^k$.
Integrating over $\Omega$ and using $m(\Omega)\le1$ gives
$\eta_n\le Ct_n^{-2k}\|\Jac(F^n)\|_\infty\,\omega^k$.
Since $t_n=t_0e^{-4\varepsilon n}$ and
$\|\Jac(F^n)\|_\infty\le\|\Jac F\|_\infty^n$, we obtain
\begin{equation}
\label{eq:approximant-density:revision}
 \eta_n\le Ce^{b_0n}\omega^k,
 \qquad b_0:=8k\varepsilon+\log\max\{1,\|\Jac F\|_\infty\}>0.
\end{equation}

\smallskip
\noindent\emph{Step 1.4. Mass bounds for the
partition elements.}
We finally estimate the $m$-mass of each set
$\Omega_{n,j,g}$ in this partition.
Fix $n\ge N$, a
non-empty element $\Omega_{n,j,g}$ of the partition,
and
$\xi\in\Omega_{n,j,g}$. For any $\zeta\in\Omega_{n,j,g}$, set
$y_\ell:=F^\ell(p(\zeta))$.
Since $y_0=g(y_n)$ and $\dist(y_n,x_n)\le2C_*R_n$,
assumption \textup{(i)} gives
$\dist(y_\ell,x_\ell)\le2C_*Aa e^{-\varepsilon n}<\delta$
for $0\le\ell<n$, by the choice of $a$.
Thus
$p(\Omega_{n,j,g})\subset B_n(x_0(\xi),\delta)$.
Consequently,
assumption \textup{(ii)} gives
\begin{equation}\label{eq:class-mass:revision}
 m(\Omega_{n,j,g})
 \le m\bigl(p^{-1}(B_n(x_0(\xi),\delta))\bigr)
 =\eta(B_n(x_0(\xi),\delta))
 \le C_Fe^{-sn}.
\end{equation}

\smallskip
\noindent{\bf Step 2. 
Comparing consecutive
regularizations at fixed parameters.}
For $u\in\Ucal$, set 
$$
I_n(\xi,u):=\int u\,dQ_{\xi,n}
\qquad
\text{and} \qquad
J_n(u):=\int u\,d\eta_n=\int_\Omega I_n(\xi,u)\,dm(\xi).$$
The density bounds for $Q_{\xi,n}$ and $\eta_n$, together
with the uniform $L^1(\omega^k)$ bound for $\Ucal$, imply that
$I_n(\xi,u)$ and $J_n(u)$ are finite.

In this step, we estimate
$|I_{n+1}(\xi,u)-I_n(\xi,u)|$
for each $\xi\in\Omega$, by expressing the two
terms as integrals over a
common domain and applying
Lemma \ref{lem:comparison}.

\smallskip
\noindent\emph{Step 2.1. An expression for
$I_{n+1}(\xi,u)-I_n(\xi,u)$.}
Fix $n\ge N$, a non-empty set $\Omega_{n,j,g}$,
$\xi\in\Omega_{n,j,g}$, and $u\in\Ucal$.
We express $I_{n+1}(\xi,u)$ and $I_n(\xi,u)$ as averages of
$u\circ g$ on $U_{n,j}$.
We already have $Q_{\xi,n}=g_*\kappa_{x_n,t_n}$.
To obtain the corresponding expression for $Q_{\xi,n+1}$,
set $b_{\xi,n}:=F^n\circ g_{\xi,n+1}$ on
$B(x_{n+1},r_0e^{-\varepsilon(n+1)})$.
Then $F\circ b_{\xi,n}=\operatorname{id}$ and
$b_{\xi,n}(x_{n+1})=x_n$, so $b_{\xi,n}$ is an inverse branch of $F$.
By Step 1.2, $\overline{K_{\xi,n+1}}$
is contained in the domain of $b_{\xi,n}$.

For $y\in\overline{K_{\xi,n+1}}$, hypothesis \textup{(i)} at
intermediate time $n$ gives
$$
 \dist(b_{\xi,n}(y),x_n)
 \le AC_*e^{\varepsilon(n+1)}e^{-\lambda}t_{n+1}
 \le AC_* t_0 e^{-3\varepsilon n}<R_n/(2C_*).
$$
Since $\dist(x_n,\partial U_{n,j}^{(2)})\ge R_n/C_*$,
we obtain $b_{\xi,n}(\overline{K_{\xi,n+1}})\subset U_{n,j}^{(2)}$.
In particular, $g\circ b_{\xi,n}$ is defined on $K_{\xi,n+1}$.
Both $g\circ b_{\xi,n}$ and $g_{\xi,n+1}$ are inverse branches of
$F^{n+1}$ there and send $x_{n+1}$ to $x_0$. Uniqueness of inverse
germs and connectedness of $K_{\xi,n+1}$ therefore give
$g_{\xi,n+1}=g\circ b_{\xi,n}$ on $K_{\xi,n+1}$.

Set $\sigma_{\xi,n}:=\kappa_{x_n,t_n}$ and
$\tau_{\xi,n}:=(b_{\xi,n})_*\kappa_{x_{n+1},t_{n+1}}$.
The measures $\sigma_{\xi,n}$ and $\tau_{\xi,n}$ are supported in
$U_{n,j}^{(2)}$, and the identity
$g_{\xi,n+1}=g\circ b_{\xi,n}$ gives
$Q_{\xi,n}=g_*\sigma_{\xi,n}$ and
$Q_{\xi,n+1}=g_*\tau_{\xi,n}$. Hence, we have
$$I_{n+1}(\xi,u)-I_n(\xi,u)
=\int_{U_{n,j}}(u\circ g)\,d(\tau_{\xi,n}-\sigma_{\xi,n}).$$

\smallskip
\noindent\emph{Step 2.2. Density bounds
of $\sigma_{\xi,n}$ and $\tau_{\xi,n}$.}
We first bound the densities of
$\sigma_{\xi,n}$ and $\tau_{\xi,n}$ with respect to $\omega^k$.
By \eqref{eq:kernel-properties}, we have
$d\sigma_{\xi,n}/d\omega^k\le Ct_n^{-2k}$.
Arguing as in \eqref{eq:inverse-branch-density},
 with $b_{\xi,n}$ in place of $g_{\xi,n}$,
we obtain
$d\tau_{\xi,n}/d\omega^k
\le Ct_{n+1}^{-2k}\|\Jac F\|_\infty$.

To
compute the densities of $\sigma_{\xi,n}$ and
$\tau_{\xi,n}$ in the coordinates of $U_{n,j}$, write
$\Phi:=\Phi_{c_{n,j}}|_{B(0,4R_n)}$ and
$J_\Phi:=d(\Phi^*\omega^k)/dV$, and set
$\widetilde\sigma_{\xi,n}:=(\Phi^{-1})_*\sigma_{\xi,n}$ and
$\widetilde\tau_{\xi,n}:=(\Phi^{-1})_*\tau_{\xi,n}$.
The uniform volume comparisons give $J_\Phi\le C_*$.
Thus, for $dV$-almost every $z\in B(0,4R_n)$, we have
$$
 \begin{aligned}
 \frac{d\widetilde\sigma_{\xi,n}}{dV}(z)
 &=\frac{d\sigma_{\xi,n}}{d\omega^k}(\Phi(z))J_\Phi(z)
 \le Ct_n^{-2k}\le Ct_{n+1}^{-2k},\\
 \frac{d\widetilde\tau_{\xi,n}}{dV}(z)
 &=\frac{d\tau_{\xi,n}}{d\omega^k}(\Phi(z))J_\Phi(z)
 \le C\|\Jac F\|_\infty t_{n+1}^{-2k}.
 \end{aligned}
$$
Absorbing the 
factor
$\max\{1,\|\Jac F\|_\infty\}$ into the constant 
$C$, the densities of
$\widetilde\sigma_{\xi,n}$ and $\widetilde\tau_{\xi,n}$ with
respect to $dV$ are bounded by $Ct_{n+1}^{-2k}$,
with $C$ independent of $n$, $j$, and $\xi$.

\smallskip
\noindent\emph{Step 2.3. Pluriharmonic cancellations.}
We next verify that the integrals of pluriharmonic functions
against $\sigma_{\xi,n}$ and $\tau_{\xi,n}$ agree.
Let $H$ be a real pluriharmonic function
on $U_{n,j}^{(3)}$.
Since $b_{\xi,n}(\overline{K_{\xi,n+1}})\subset U_{n,j}^{(2)}$,
the function $H\circ b_{\xi,n}$ is pluriharmonic on a
neighbourhood of $\overline{K_{\xi,n+1}}$.
The radial mean value identity \eqref{eq:kernel-properties} gives
\begin{equation}\label{eq:average-cancellation:revision}
 \int H\,d\tau_{\xi,n}=H(x_n)=\int H\,d\sigma_{\xi,n},
\end{equation}
which proves the claim.

\smallskip
\noindent\emph{Step 2.4. Current masses.}
We continue to denote 
$\Phi:=\Phi_{c_{n,j}}|_{B(0,4R_n)}$.
Set $T_u:=\omega+\ddc u$ and
$\omega_{n,j}:=(\Phi^{-1})^*\omega_{\mathrm{eucl}}$, and define
\begin{equation}\label{eq:local-current-mass}
 M_{n,j,g}(u):=\int_{U_{n,j}}g^*(T_u+\omega)
 \wedge\omega_{n,j}^{k-1}.
\end{equation}
Since $g$ is biholomorphic onto its image, $u\circ g\circ\Phi$
is locally integrable on $B(0,4R_n)$ and
$\ddc(u\circ g\circ\Phi)=\Phi^*g^*T_u-\Phi^*g^*\omega$.
Since $g$ is holomorphic on a neighbourhood of
$\overline{U_{n,j}}$, the positive currents $\Phi^*g^*T_u$ and
$\Phi^*g^*\omega$ have finite mass on $B(0,4R_n)$, with
$$
 \|\Phi^*g^*T_u\|_{B(0,4R_n)}
 +\|\Phi^*g^*\omega\|_{B(0,4R_n)}=M_{n,j,g}(u).
$$

\smallskip
\noindent\emph{Step 2.5. Application of Lemma \ref{lem:comparison} and conclusion.}
Steps
2.2--2.4 
show that the measures $\widetilde\sigma_{\xi,n}$ and
$\widetilde\tau_{\xi,n}$ are supported in $B(0,2R_n)$, have
densities at most $Ct_{n+1}^{-2k}$ with respect to $dV$, and
have equal integrals against every real
pluriharmonic function
$\widetilde H$ on $B(0,3R_n)$.
The last property
follows from
\eqref{eq:average-cancellation:revision} applied to
$H=\widetilde H\circ\Phi^{-1}$.
Moreover, we have
$\ddc(u\circ g\circ\Phi)
 =\Phi^*g^*T_u-\Phi^*g^*\omega$,
and the sum of the masses of these two positive closed currents
on $B(0,4R_n)$ is $M_{n,j,g}(u)$.
We may therefore apply Lemma \ref{lem:comparison}
with $r=R_n$ and $M=Ct_{n+1}^{-2k}$.
This gives
\begin{equation}
\label{eq:pointwise-comparison:revision}
 |I_{n+1}(\xi,u)-I_n(\xi,u)|
 \le CR_n^2t_{n+1}^{-2k}M_{n,j,g}(u)
 \le Ce^{(8k-4)\varepsilon n}M_{n,j,g}(u).
\end{equation}
Observe that the quantity $M_{n,j,g}(u)$ depends on the fixed branch $g$,
not on the choice of $\xi\in\Omega_{n,j,g}$.

\smallskip
\noindent{\bf Step 3. A bound for
$|J_{n+1}(u)-J_n(u)|$.}
In this step we first
prove that
the quantities
$M_{n,j,g}(u)$ defined in \eqref{eq:local-current-mass}
satisfy
\begin{equation}\label{eq:total-current-mass:revision}
 \sum_j\sum_{g\in\mathcal G_{n,j}}M_{n,j,g}(u)
 \le Cd^{(k-1)n}
\end{equation}
uniformly for $u\in\Ucal$, and use this bound to estimate $|J_{n+1}(u)-J_n(u)|$.

\smallskip
\noindent\emph{Step 3.1. The current-mass bound \eqref{eq:total-current-mass:revision}.}
For fixed $n,j$, distinct branches
$g,h\in\mathcal G_{n,j}$ have disjoint images. Indeed, if
$g(y)=h(z)$ for some $y,z\in U_{n,j}$, applying $F^n$ gives
$y=z$. Thus $g(y)=h(y)$, and uniqueness of inverse germs and
connectedness of $U_{n,j}$ imply $g=h$ on $U_{n,j}$.

For every global positive closed $(1,1)$-current $S$, we have
\begin{equation}
\label{eq:current-domination:revision}
 \sum_{g\in\mathcal G_{n,j}}g^*S
 \le \left.(F^n)_*S\right|_{U_{n,j}}.
\end{equation}
Indeed, let $\Psi$ be a smooth positive $(k-1,k-1)$-form compactly
supported in $U_{n,j}$, extended by zero to $\PP^k$.
Since $\Psi$ is compactly supported in $U_{n,j}$, each
$g_*\Psi$ has compact support in $g(U_{n,j})$ and extends smoothly
by zero to $\PP^k$. On $g(U_{n,j})$, we have
$g_*\Psi=(F^n)^*\Psi$. Since the images of distinct branches are
disjoint, it follows that
$\sum_g g_*\Psi\le(F^n)^*\Psi$.
Therefore,
$$
 \sum_g\langle g^*S,\Psi\rangle
 =\Big\langle S,\sum_g g_*\Psi\Big\rangle
 \le\langle S,(F^n)^*\Psi\rangle
 =\langle(F^n)_*S,\Psi\rangle,
$$
which gives \eqref{eq:current-domination:revision}. Observe that no smoothness of $S$ is assumed.

Set $S:=T_u+\omega$. Since $\omega_{n,j}$ is uniformly comparable
to $\omega$, \eqref{eq:current-domination:revision} and the uniformly
bounded overlap of the domains $U_{n,j}$ give
$$
\begin{aligned}
 \sum_j\sum_{g\in\mathcal G_{n,j}}M_{n,j,g}(u)
 &\le
 C\sum_j\int_{U_{n,j}}
 \left(\sum_{g\in\mathcal G_{n,j}}g^*S\right)
 \wedge\omega^{k-1} \\
 &\le
 C\sum_j\int_{U_{n,j}}(F^n)_*S\wedge\omega^{k-1} 
 \le
 C\int_{\PP^k}(F^n)_*S\wedge\omega^{k-1}
 =2Cd^{(k-1)n}.
\end{aligned}$$
The equality follows from the projection formula and cohomology,
using $[T_u+\omega]=2[\omega]$,
$(F^n)^*[\omega^{k-1}]
=d^{(k-1)n}[\omega^{k-1}]$,
and
$\int_{\PP^k}\omega^k=1$.
This proves \eqref{eq:total-current-mass:revision}.
Observe
that
the sum is controlled by current
growth, with no factor for the number of branches or domains $U_{n,j}$.

\smallskip
\noindent\emph{Step 3.2. Estimating $|J_{n+1}(u)-J_n(u)|$.}
We now combine \eqref{eq:total-current-mass:revision} with
\eqref{eq:class-mass:revision}
to bound $|J_{n+1}(u)-J_n(u)|$.
Set
$a_0:=s-(k-1)\log d-(8k-4)\varepsilon>0$.
Integrating \eqref{eq:pointwise-comparison:revision} over each element
$\Omega_{n,j,g}$ of the partition and summing, then using
\eqref{eq:class-mass:revision} and \eqref{eq:total-current-mass:revision}, yields,
for every $n\ge N$,
\begin{equation}\label{eq:summable-comparison:revision}
 |J_{n+1}(u)-J_n(u)|
 \le Ce^{(8k-4)\varepsilon n}
 \sum_{j,g}m(\Omega_{n,j,g})M_{n,j,g}(u)
 \le Ce^{-a_0n}.
\end{equation}

\smallskip
\noindent{\bf Step 4. H\"older regularity.}
Together with the equality of the masses and the weak convergence
established in Step 1, the density bound
\eqref{eq:approximant-density:revision} and the comparison estimate
\eqref{eq:summable-comparison:revision} show that the hypotheses
of Lemma \ref{lem:holder-limit} hold
with $a=a_0$ and $b=b_0$. Therefore $\eta$ integrates every
quasi-plurisubharmonic function and satisfies
\eqref{eq:regularization-holder:revision} with
$$
 \beta=\frac{a_0}{a_0+b_0}\in(0,1).
$$
Since $a_0=s-(k-1)\log d-(8k-4)\varepsilon$
 and $b_0=8k\varepsilon+
 \log\max\{1,\|\Jac F\|_\infty\}$,
the exponent has the desired dependence. The final statement
about the super-potential follows from
\eqref{eq:regularization-holder:revision} under the identification
recalled after Theorem \ref{thm:submeasure}.
\end{proof}

\section{Applications of the regularization criterion}
\label{sec:applications}

In this section,
we apply Theorem \ref{prop:regularization} to prove Theorems
\ref{thm:submeasure} and \ref{thm:bowen}, and deduce Corollaries
\ref{cor:entropy} and \ref{cor:moderate}.
We also give a direct proof in dimension one and an example
showing that global quasi-plurisubharmonic integrability can fail.

\subsection{Properties of large entropy measures}\label{subsec:restrictions}
\begin{proof}[Proof of Theorem \ref{thm:submeasure}]
We first construct a submeasure 
of $\mu$
to which Theorem
\ref{prop:regularization} applies. Let $\chi_{\min}>0$ be the smallest
Lyapunov exponent of $\mu$; see Section \ref{subsec:branches}.
Choose $(k-1)\log d<s<h_\mu(F)$ and
$0<\varepsilon<\chi_{\min}/2$ with
$(8k-4)\varepsilon<s-(k-1)\log d$.
Fix $0<\theta<1$, apply Lemma \ref{lem:good-set} with mass error
$\theta/4$, and let $\widehat E$
be the resulting Borel set. Set
$\Omega:=\widehat E$, $m:=\widehat\mu|_{\widehat E}$,
$p:=\pi_0|_{\widehat E}$, and
$\eta:=p_*m$. Then $\eta\le\mu$ and $\eta(\PP^k)>1-\theta/4$.
Lemma \ref{lem:good-set} \textup{(i)} gives the branch estimates,
and \textup{(ii)} gives
$\eta(B_n(x_0,\delta))\le\mu(B_n(x_0,\delta))\le e^{-sn}$.
Theorem \ref{prop:regularization} gives integrability of
quasi-plurisubharmonic functions and the absolute H\"older estimate
\eqref{eq:regularization-holder:revision} for $\eta$. Since $s$ and $\varepsilon$
were chosen independently of $\theta$, so is the exponent $\beta$.

We now choose a compact set $E_\theta$ with
$\mu(E_\theta)>1-\theta$ such that $\one_{E_\theta}\mu\le2\eta$.
Choose a Borel density $f$ with $\eta=f\mu$ and $0\le f\le1$, and denote
$A_\theta:=\{f\ge1/2\}$. Since $1-f>1/2$ on $A_\theta^c$,
we have $\mu(A_\theta^c)\le2\int_{\PP^k}(1-f)\,d\mu<\theta/2$.
Since $\mu(A_\theta)>1-\theta/2$, inner regularity gives a
compact set $E_\theta\subset A_\theta$ with
$\mu(E_\theta)>1-\theta$.
The domination $\eta_\theta:=\one_{E_\theta}\mu\le2\eta$
gives integrability and
$\int|u-v|\,d\eta_\theta\le2\int|u-v|\,d\eta$.
Thus the absolute H\"older estimate holds for every
$\eta_\theta$
with the same exponent $\beta$.
\end{proof}

\begin{remark}\label{rem:holder-exponent}
Set $\Gamma_F:=\log\max\{1,\|\Jac F\|_\infty\}$. In
Theorem \ref{thm:submeasure}, every exponent
$$
 0<\beta<
 \frac{h_\mu(F)-(k-1)\log d}
{h_\mu(F)-(k-1)\log d+\Gamma_F}
$$
can be used for all $0<\theta<1$, with a possibly
$\theta$-dependent constant $C_\theta$.
Indeed, Step 4 of the proof of Theorem \ref{prop:regularization}
gives the exponent
$\beta_0=(s-(k-1)\log d-(8k-4)\varepsilon)/
(s-(k-1)\log d+4\varepsilon+\Gamma_F)$.
Choose $s<h_\mu(F)$ sufficiently close to $h_\mu(F)$ and then
$\varepsilon>0$ sufficiently small, independently of $\theta$,
so that the conditions on $s$ and $\varepsilon$ in the proof of
Theorem \ref{thm:submeasure} hold and $\beta_0>\beta$.
The estimate with exponent $\beta_0$ implies the estimate with
exponent $\beta$, since
$\max\{t,t^{\beta_0}\}\le\max\{t,t^\beta\}$ for $t\ge0$.
\end{remark}

\begin{proof}[Proof of Corollary \ref{cor:entropy}]
We prove \textup{(i)} using integrability
of the restrictions
of $\mu$ to sets of measure
arbitrarily close to one.
Let $P$ be a Borel pluripolar set.
Choose $u\in\Ucal$ with $P\subset\{u=-\infty\}$.
Integrability of $u$ with respect to
the restrictions in Theorem \ref{thm:submeasure} gives
$\eta_\theta(P)=0$ and hence
$\mu(P)\le(\mu-\eta_\theta)(\PP^k)<\theta$ for every $0<\theta<1$.
Letting $\theta\to0$ gives $\mu(P)=0$.

We deduce \textup{(ii)} by ergodic decomposition and the variational
principle. If $\nu(P)=1$, then almost every component $\rho$
in the ergodic decomposition of $\nu$ satisfies $\rho(P)=1$,
so \textup{(i)} gives
$h_\rho(F)\le(k-1)\log d$. Integrating the component entropies gives
$h_\nu(F)\le(k-1)\log d$. If $K$ is non-empty, compact, forward
invariant, and pluripolar, every invariant probability measure for
$F|_K$ satisfies this entropy bound. The variational
principle therefore gives $\htop(F|_K)\le(k-1)\log d$.
\end{proof}

\begin{proof}[Proof of Corollary \ref{cor:moderate}]
We prove exponential integrability with an exponent independent of
$\theta$, using the truncation argument in
\cite[proof of Proposition 3.9]{DN14}.
Skoda's integrability theorem
\cite{Sko72},
in the uniform form of \cite[Theorem A.22]{DS10},
gives constants
$c,C>0$
such that
$\|\max(u,-t)-u\|_1\le Ce^{-ct}$ for $u\in\Ucal$ and $t\ge1$.
Apply \eqref{eq:holder-restriction} to
$u$ and $u_t:=\max(u,-t)$. For all sufficiently large $t$,
we have
$$
 \eta_\theta(\{u<-t-1\})
 \le\int_{\PP^k}(u_t-u)\,d\eta_\theta
 \le C_\theta e^{-c\beta t}.
$$
Integration of the sublevel bound gives
$\sup_{u\in\Ucal}\int_{\PP^k} e^{-\alpha u}\,d\eta_\theta<\infty$
for every $0<\alpha<c\beta$. Both $c$ and $\beta$ are independent
of $\theta$, so the same $\alpha$ works for all restrictions.

We obtain the Monge--Amp\`ere representation with a common
H\"older exponent by applying
\cite[Theorem 1.1]{DN14} to $\eta_\theta/\eta_\theta(\PP^k)$.
Normalization preserves the super-potential exponent $\beta$ in
\eqref{eq:holder-restriction}. The exponent of the resulting
potential depends only on $k$ and $\beta$, so it can be chosen
independently of $\theta$.
\end{proof}

\subsection{Frostman measures and Bowen entropy}
\label{sec:bowen}
Since uniform expansion gives inverse branches on fixed-radius balls,
a uniform Bowen ball mass bound suffices to apply
Theorem \ref{prop:regularization}.
\begin{proposition}\label{prop:frostman}
Let $F\colon\PP^k\to\PP^k$ be a holomorphic endomorphism of
algebraic degree $d\ge2$, with $k\ge1$, and let $K\subset\PP^k$
be non-empty and compact, with $F(K)\subset K$.
Suppose that $F$ is uniformly expanding on $K$. Let $\sigma$ be a
Borel probability measure on $\PP^k$, not necessarily invariant,
with $\sigma(K)=1$. If there are
$\delta,C_F>0$, an integer $N\ge1$, and $s>(k-1)\log d$ such that
$\sigma(B_n(x,\delta))\le C_Fe^{-sn}$ for every $x\in K$ and
$n\ge N$, then $\sigma$ has a H\"older continuous super-potential.
In particular, it integrates every quasi-plurisubharmonic function
and is moderate.
\end{proposition}

\begin{proof}
We first obtain uniform inverse branches along orbits in $K$.
Uniform expansion gives an integer $m\ge1$ with
$\sup_{x\in K}\|(D(F^m)_x)^{-1}\|<1$.
Compactness of $K$ and the inverse function theorem give local
inverse branches of $F^m$ on balls of a common radius, with a
common contraction factor strictly less than one.
The inverse branches of the remaining iterates $F^r$, $0\le r<m$,
have uniformly bounded Lipschitz constants on smaller balls of a
common radius. Writing $n=qm+r$, with $q\ge0$ and $0\le r<m$,
we compose the inverse branches of $F^m$ with the remaining inverse
branches of $F^r$. Absorbing the bounds for $0\le r<m$ into a
constant gives $r_0,\lambda>0$ and $A\ge1$
such that $g_{x,n}$, the inverse branch sending $F^nx$ to $x$,
is defined on $B(F^nx,r_0)$ and satisfies
$\Lip(F^j\circ g_{x,n})\le Ae^{-\lambda(n-j)}$
for $x\in K$ and $0\le j\le n$.

We now apply Theorem \ref{prop:regularization} with $\Omega=K$,
parameter measure $\sigma|_K$, and $p$ the inclusion.
Choose $\varepsilon>0$ with $(8k-4)\varepsilon<s-(k-1)\log d$.
Restricting the branches to $B(F^nx,r_0e^{-\varepsilon n})$
gives hypothesis \textup{(i)}, and the assumed Bowen ball bound gives
\textup{(ii)}. Theorem \ref{prop:regularization} then
gives integrability and the absolute
H\"older estimate for $\sigma$, and 
\cite[Proposition 3.9]{DN14} gives moderateness.
\end{proof}

Following \cite{FH12}, we use variable length covers to measure the
orbit complexity of sets that need not be compact or invariant. For $Z\subset \PP^k$, $s\ge0$,
$\delta>0$, and $N\ge1$, denote
$$
 \mathcal M_{N,\delta}^s(Z)
 :=\inf\left\{\sum_i e^{-sn_i}:
 Z\subset\bigcup_i B_{n_i}(x_i,\delta),\ n_i\ge N\right\},
$$
where the infimum is over countable covers with $x_i\in \PP^k$.
Set $\mathcal M_\delta^s(Z):=\lim_{N\to\infty}\mathcal M_{N,\delta}^s(Z)$
and define the \emph{Bowen entropy} by
$h_{\mathrm B}(F,Z):=\lim_{\delta\to 0}
\inf\{s\ge0:\mathcal M_\delta^s(Z)=0\}$.
Bowen entropy agrees with topological entropy on compact invariant
sets and is countably stable. The use of variable lengths
distinguishes Bowen entropy from upper capacity entropy, which uses
a common orbit length; see \cite[Section 2]{FH12}.

\begin{proof}[Proof of Theorem \ref{thm:bowen}]
Suppose on the contrary that
$h_{\mathrm B}(F,P\cap K)>(k-1)\log d$.
Since $P\cap K$ is Borel and $\htop(F)<\infty$, the compact
approximation theorem \cite[Theorem 1.2 (ii)]{FH12}
gives a compact set $L\subset P\cap K$ with
$h_{\mathrm B}(F,L)>(k-1)\log d$.
Choose $(k-1)\log d<s<s+\gamma<h_{\mathrm B}(F,L)$.
For sufficiently small $\delta>0$ and all sufficiently large $N$,
we have $\mathcal M_{N,6\delta}^{s+\gamma}(L)>0$.
Applying \cite[Proposition 3.2 and Lemma 3.4]{FH12}, we obtain
a Borel probability measure $\sigma$ with $\sigma(L)=1$ and
$\sigma(B_n(x,\delta))\le C_Fe^{-sn}$ for every $x\in\PP^k$
and $n\ge N$, for some $C_F>0$.

Since $\sigma(K)=1$, Proposition \ref{prop:frostman} implies that
every $u\in\Ucal$ is $\sigma$-integrable. Choose $u\in\Ucal$ with $P\subset\{u=-\infty\}$. 
The equality $\sigma(L)=1$, with
$L\subset P$, forces $\int_{\PP^k}(-u)\,d\sigma=\infty$, which is a contradiction.
\end{proof}

\subsection{The one-dimensional case}
\label{subsec:dimension-one}
In dimension one, Theorem \ref{thm:submeasure} is an immediate
consequence of classical dimension theory and potential theory.
Let $F\colon\PP^1\to\PP^1$ be a rational map of degree at least
two and let $\mu$ be an ergodic invariant probability measure with
$h_\mu(F)>0$. By Ma\~n\'e's dimension formula
\cite[Theorem A]{Man88}, the measure $\mu$ is exact-dimensional,
with
$\dim\mu={h_\mu(F)}/{\chi_\mu}>0$,
where $\chi_\mu$ is its Lyapunov exponent.
Fix $0<a<\dim\mu$. Exact dimensionality implies that, for every
$0<\theta<1$, there are a compact set
$E_\theta\subset\PP^1$, with $\mu(E_\theta)>1-\theta$, and a
constant $C_\theta>0$ such that
$\mu(B(x,r))\le C_\theta r^a$
for every 
$x\in E_\theta$
and $r>0$. For $\eta_\theta:=\one_{E_\theta}\mu$, the same estimate, after
changing $C_\theta$, holds for balls with arbitrary centers: if
$B(y,r)$ meets $E_\theta$, choose
$x\in E_\theta\cap B(y,r)$ and use
$B(y,r)\subset B(x,2r)$.

By \cite[Proposition 3.8]{DN14}, the measure $\eta_\theta$ has a
H\"older continuous super-potential. Since the exponent in that
proposition depends only on $a$, it can be
chosen independently of $\theta$. The estimate
\eqref{eq:holder-restriction} then follows from
\cite[Lemma 3.5]{DN14}.

In dimension one, Theorem \ref{thm:bowen}
holds without
uniform expansion. Fix $L>\max\{1,\Lip(F)\}$. Polar sets have Hausdorff dimension zero
\cite[Theorem 5.9.6]{AG01}, and comparison of ordinary ball covers
with Bowen ball covers gives
$h_{\mathrm B}(F,P)\le(\log L)\dim_{\mathrm H}P=0$
for every Borel polar set $P\subset\PP^1$, where
$\dim_{\mathrm H}$ denotes Hausdorff dimension.

\subsection{Failure of global integrability}
The following standard renewal example, included for completeness,
shows that global integrability can fail already in dimension one,
even for mixing measures with entropy arbitrarily close to the
maximum.

\begin{proposition}
\label{prop:nonintegrable}
Let $Q\colon\PP^1\to\PP^1$ be given by
$Q[Z_0:Z_1]=[Z_0^2:Z_1^2]$
and set
$$
 u([Z_0:Z_1])
 :=
 \log\frac{|Z_1-Z_0|}
 {\sqrt2\bigl(|Z_0|^2+|Z_1|^2\bigr)^{1/2}}.
$$
Then $u\in\Ucal$.
For
every $\varepsilon>0$, there is a mixing $Q$-invariant probability
measure $\nu$ supported on the uniformly expanding repeller
$\mathbb S^1:=\{w\in\C:|w|=1\}$, where $w=Z_1/Z_0$,
satisfying
$$
 h_\nu(Q)>\log2-\varepsilon,
 \qquad
 \nu\{u=-\infty\}=0,
 \qquad
\text{and}
\qquad
 \int_{\PP^1}(-u)\,d\nu=\infty.
$$
\end{proposition}

\begin{proof}
The statement follows from a standard stationary-renewal
construction; see \cite[Chapter XI]{Fel71} for the equilibrium
renewal process and \cite{Abr59} for the entropy formula for induced
transformations.
Denote $\mathbb N_0:=\{0,1,2,\ldots\}$.
For $0<\delta<1$, let $q=(q_\ell)_{\ell\ge0}$ be the probability
distribution 
on $\mathbb N_0$
defined by
$q_\ell:=(1-\delta)2^{-(\ell+1)}+\delta c(\ell+1)^{-3}$, where
$c^{-1}:=\sum_{\ell\ge0}(\ell+1)^{-3}$.
Set
$$
 \overline L:=\sum_{\ell\ge0}(\ell+1)q_\ell.
$$
Let $(X_j)_{j\ge0}$ be
the stationary Markov chain on $\mathbb N_0$ with deterministic
transitions $j\mapsto j-1$ for $j\ge1$ and transition probabilities
$\Pr(0\mapsto\ell)=q_\ell$.
Let $\varrho_\delta$ be the law on $\{0,1\}^{\mathbb N_0}$ of
$a_j:=\one_{\{X_j=0\}}$. 
This is a stationary binary renewal process in which
the numbers of zeros between successive $1$'s are
independent and have distribution $q$.
This chain is irreducible, aperiodic, and positive recurrent; hence
it is mixing by \cite[Theorem 1.8.3]{Nor97}. Therefore
$\varrho_\delta$ is mixing.
For $a=(a_j)_{j\ge0}$, set
$R(a):=\min\{j\ge0:a_j=1\}$, with $R(a)=\infty$ if no such
index exists. Standard renewal theory gives
$$
 \varrho_\delta
 (\{R=j\})
 =\frac1{\overline L}\sum_{\ell\ge j}q_\ell,
 \qquad \text{hence}\qquad
 \mathbb E_{\varrho_\delta}R
 =\frac1{2\overline L}
 \sum_{\ell\ge0}\ell(\ell+1)q_\ell
 =\infty.
$$
On the other hand, Abramov's formula \cite{Abr59}
gives
$h_{\varrho_\delta}(\sigma)={H(q)}/{\overline L}$,
where $\sigma$ is the left shift on $\{0,1\}^{\mathbb N_0}$ and
$H(q):=-\sum_{\ell\ge0}q_\ell\log q_\ell$.
As $\delta\to0$, we have $ \overline L\to2$. Moreover,
the geometric distribution has entropy $2\log2$,
and the polynomial-tail distribution has finite
entropy. The entropy bounds for convex mixtures
therefore give $H(q)\to2\log2$.
Hence we have
$h_{\varrho_\delta}(\sigma)\to\log2$.
Let
$$
 \Pi(a)=
 \exp\left(2\pi i\sum_{j\ge0}a_j2^{-(j+1)}\right)
$$
be the binary coding and set
$\nu_\delta:=\Pi_*\varrho_\delta$.
Then $\nu_\delta$ is mixing for $w\mapsto w^2$. The coding is
one-to-one outside the dyadic ambiguities, which have
$\varrho_\delta$-measure zero, and therefore we have
$h_{\nu_\delta}(w\mapsto w^2)
 =h_{\varrho_\delta}(\sigma)$.
Moreover, $\nu_\delta (\{1\})=0$, while
$|\Pi(a)-1|\le 2\pi 2^{-R(a)}$
implies
$$
 \int_{\mathbb S^1}
 \log\frac{2}{|w-1|}\,d\nu_\delta(w)=\infty.
$$

By the Poincar\'e--Lelong formula, we have
$\ddc u=\delta_1-\omega$, where $\delta_1$ is the Dirac measure
at the point $w=1$, so
$u\in\PSH(\PP^1,\omega)$ and
$\{u=-\infty\}=\{1\}$.
Moreover, $u\le0$ by the Cauchy--Schwarz inequality, with
$u([1:-1])=0$. Hence $u\in\Ucal$.
On $\mathbb S^1$ we have
$-u(w)=\log ({2} /|w-1|)$.
Choosing $\delta>0$ sufficiently small and setting
$\nu:=\nu_\delta$ proves all the assertions.
\end{proof}

\section{Positive-entropy measures in saddle two-dimensional settings}
\label{sec:surface}
The following result is independent of the regularization criterion
of Theorem \ref{prop:regularization}. Its proof combines standard
consequences of the Ledrappier--Young
theory \cite{LY85}
with a pluripolarity
criterion of Dujardin \cite{Duj05}.

\begin{theorem}
\label{thm:surface}
Let $F$ be in one of the following settings:
\begin{enumerate}[label=\textup{(\roman*)},leftmargin=2em]
\item
$F\colon X\to X$ is a holomorphic automorphism
of a connected compact K\"ahler surface $X$;
\item
$F\colon\C^2\to\C^2$ is a polynomial automorphism
of H\'enon type;
\item
$F$ is an invertible horizontal-like map on a bidisk
$D\subset\C^2$ \cite{Duj04HL}.
\end{enumerate}
In case \textup{(iii)}, let $K_F\Subset D$ denote the compact
maximal invariant set, consisting of the points whose full orbits
are defined and remain in $D$.
Write $Y:=X$, $\C^2$, or $D$, respectively.
Let $\mu$ be an ergodic $F$-invariant probability measure with
$h_\mu(F)>0$, assumed to be compactly supported in case
\textup{(ii)} and supported on $K_F$ in case \textup{(iii)}.
Then $\mu$ gives zero mass to every Borel pluripolar set in $Y$.
\end{theorem}

For the canonical measures of maximal entropy in these
settings, the conclusion follows from their representation
as wedge products of Green currents with continuous
local potentials;
see
\cite[Section 5.2]{Can14} for compact K\"ahler surfaces,
\cite{BS91,BLS93} for polynomial automorphisms
of H\'enon type, and \cite{Duj04HL} for
invertible horizontal-like maps.

\begin{proof}
The argument is local near $\supp(\mu)$ and applies in all three
settings. In each case, $F$ and $F^{-1}$ are holomorphic on
neighbourhoods of the compact support, so their derivatives are
bounded there. Thus the integrability hypotheses needed for
Ruelle's inequality and Pesin theory hold.

Applying Ruelle's entropy inequality to $F$ and $F^{-1}$ shows that
$\mu$ has Lyapunov exponents
$\chi^u>0>\chi^s$,
each of real multiplicity two. By Pesin theory, the local stable and
unstable manifolds through $\mu$-almost every point are complex
disks. Moreover, a set of full measure can be covered by countably
many compact Pesin boxes in which the local stable and unstable
plaques have uniform size, depend continuously on the base point,
and form two measurable families of pairwise disjoint holomorphic
graphs. We refer to
\cite[Sections 2.4--2.6 and Section 4]{BLS93}, \cite[Section 5]{Duj04HL},
and
\cite[Sections 8.1.2--8.1.3]{Can14}.

We first recall a standard consequence of the Ledrappier--Young theory.
Applied to $F$ and $F^{-1}$, its conditional dimension theorem gives
positive pointwise dimension to the conditional measures of $\mu$
along both unstable and stable Pesin manifolds; more precisely, the
corresponding dimensions are
 $h_\mu(F) /\chi^u$ and 
 $h_\mu(F) / |\chi^s|$,
respectively; see
\cite[Proposition 7.3.1 and Theorem C']{LY85}.
Since polar subsets of a complex disk have Hausdorff dimension zero
\cite[Theorem 5.9.6]{AG01}, these leafwise conditional measures give
zero mass to polar sets. The same remains true after restricting
$\mu$ to a Pesin box: the resulting plaque conditionals are
normalized restrictions of the corresponding leafwise conditional
measures; see
\cite[Section 2.4, Proposition 2.5 and the preceding discussion]{BLS93}.

Suppose now that $P\subset Y$
is Borel and pluripolar and that
$\mu(P)>0$. Choose a connected coordinate neighborhood $U$ and a
plurisubharmonic function $\varphi\not\equiv-\infty$ on $U$ such
that
 $\mu(P\cap U)>0$ and
 $P\cap U\subset\{\varphi=-\infty\}$.
We can then choose a compact Pesin box $R\subset U$ such that
$\mu(P\cap R)>0$ and all its local stable and unstable plaques are
contained in $U$. Set
$\nu:={\one_R\mu} /{\mu(R)}$,
 $E:=P\cap R$,
 and write $W_R^u(x)$ and $W_R^s(x)$ for the local holomorphic
plaques through $x$. Denote by $\nu_x^u$ and $\nu_x^s$ the conditional
probabilities of $\nu$ on the partitions of $R$ by the intersections
with these plaques.

Disintegrating $\nu$ along unstable plaques gives a set of positive
quotient measure on which $\nu_x^u(E)>0$. Choose Borel versions
of the conditional probabilities that are constant on each
partition element. Let $N^u\subset R$ be a saturated Borel set
with $\nu(N^u)=0$ such that, outside $N^u$, the conditional
probabilities vanish on polar sets. Then
$ B:=\{x\in R\setminus N^u:\nu_x^u(E)>0\}$
is measurable and has positive $\nu$-measure.
For every plaque
selected by $B$, the set $E\cap W_R^u(x)$ is nonpolar. Since the
restriction of $\varphi$ to $W_R^u(x)$ is subharmonic or
identically $-\infty$ and has $E\cap W_R^u(x)$ in its pole set,
this restriction is identically $-\infty$. Thus every full
unstable plaque selected by $B$ lies in
$\{\varphi=-\infty\}$.

We now disintegrate $\nu$ along stable plaques. There is a stable
plaque $D$ whose conditional measure gives positive mass to $B$ and
vanishes on polar subsets of $D$. 
Cover $B\cap D$ by countably many smaller flow boxes in which $D$
is transverse to the unstable plaques and these plaques are graphs
over a common disk. At least one of these boxes meets $B\cap D$
in a set of positive stable conditional measure. By inner
regularity, there is a compact set $A\subset B\cap D$, contained
in this box, with positive stable conditional measure. In
particular, $A$ is nonpolar in the complex disk $D$.

The unstable plaques through the points of $A$ form a continuous
family of pairwise disjoint holomorphic graphs, and their union is
contained in $\{\varphi=-\infty\}$. After straightening the
transversal $D$, \cite[Lemma 6.5]{Duj05} implies that the parameter
set $A$ is polar. This contradicts the preceding paragraph.
Therefore $\mu(P)=0$.
\end{proof}

By ergodic decomposition, Theorem \ref{thm:surface} also implies
that every $F$-invariant probability measure satisfying the same
support assumptions and giving full mass to a Borel pluripolar set
has zero entropy. By the variational principle, every non-empty
compact pluripolar set $K\subset Y$ satisfying $F(K)\subset K$
has
$h_{\mathrm{top}}(F|_K)=0$.
In case \textup{(iii)}, we assume in addition that $K\subset K_F$.

\begin{remark}
\label{rem:higher-dimensional-saddle}
The proof of Theorem \ref{thm:surface} uses that the stable and
unstable Pesin manifolds are complex curves, on which polar sets
have Hausdorff dimension zero. In higher dimension, positive
conditional dimension alone does not exclude mass on pluripolar
subsets of the Pesin manifolds, so the argument does not
directly extend.
\end{remark}

\appendix

\section{Compact K\"ahler manifolds and
polynomial-like maps}
\label{sec:kahler}

Theorems \ref{thm:submeasure} and \ref{thm:bowen} and
Corollaries \ref{cor:entropy} and \ref{cor:moderate}
extend to the following settings, with entropy threshold
$\log\Lambda$ in place of $(k-1)\log d$.

\smallskip
\noindent{\bf Compact K\"ahler manifolds.}
Let $F\colon X\to X$ be a surjective holomorphic self-map of a
connected compact K\"ahler manifold $(X,\omega)$ of dimension
$k\ge1$, normalized by $\int_X\omega^k=1$.
Denote by $\lambda_j(F)$ the dynamical degrees and assume
$\lambda_k(F)>\Lambda:=\lambda_{k-1}(F)$.
The conclusions hold on $X$, with the positivity conditions
and norms defined using $\omega$.
For super-potentials on compact K\"ahler
manifolds, see \cite{DS10Aut}.
In dimension two, the remaining case
$\lambda_2(F)\le\lambda_1(F)$ is covered by Gromov's entropy
bound \cite{Gro03}. Thus, for every surjective holomorphic
self-map of a compact K\"ahler surface, an invariant measure
supported on a compact pluripolar set has entropy at most
$\log\lambda_1(F)$, answering a question in
\cite[Section 2]{dTV15}.

\smallskip
\noindent{\bf Polynomial-like maps.}
Let $f\colon U\to V$ be a proper holomorphic map, where
$U\Subset V\subset\C^k$ and $V$ is convex, and let
$K_f:=\bigcap_{n\ge1}f^{-n}(U)$ be its filled Julia set.
Assume that the topological degree satisfies
$d_t>\Lambda:=d_{k-1}^*(f)$, where $d_{k-1}^*(f)$ is the
$*$-dynamical degree controlling the growth of pushforwards
of positive closed $(1,1)$-currents; see \cite{DS03,DS10}.
The conclusions hold for measures supported on $K_f$,
viewed as measures on $\PP^k$ through the standard affine chart.
Pluripolar sets are understood in $V$, and the expanding compact
set in Theorem \ref{thm:bowen} is assumed to lie in $K_f$.

\smallskip

We explain the changes needed in 
the proof of Theorem \ref{prop:regularization} 
in both cases.
The chart, covering, radial averaging, and local comparison
arguments apply on a finite atlas of $X$, or on a fixed
neighbourhood $K_f\Subset W\Subset U$.
In the polynomial-like case, the transported kernels and their
 forward images remain in $W$, where the Jacobian
is uniformly bounded.

\smallskip

The main change in both settings concerns the global
current-mass estimate
\eqref{eq:total-current-mass:revision}.
On $X$, cohomology gives
$\|(F^n)_*(T_u + \omega)\|\le C_\delta\delta^n$
for every $\delta>\Lambda$, uniformly in
$u\in\PSH(X,\omega)$.
In the polynomial-like case, the definition of $d_{k-1}^*(f)$
gives $\|(f^n)_*S\|_W\le C_\delta\delta^n$
for every $\delta>\Lambda$ and positive closed
$(1,1)$-currents $S$ of uniformly bounded mass on $U$.
By the same current domination argument in Step 3.1
of the proof of Theorem \ref{prop:regularization},
these bounds yield
\[
\sum_j\sum_{g\in\mathcal G_{n,j}}M_{n,j,g}(u)
\le C_\delta\delta^n,
\qquad \delta>\Lambda,
\]
uniformly in $u$, which replaces
\eqref{eq:total-current-mass:revision}.

\smallskip

We next verify the inverse branch hypotheses for invariant
measures of entropy greater than $\log\Lambda$.
In the compact K\"ahler case, the proof of the
entropy--Lyapunov estimate in
\cite[Section 6]{Dup12} still applies with minimal modifications:
the backward graph transform is local, and
\cite[Lemma 5]{dT08} gives the required
mixed pullback volume estimates.
By log-concavity of the dynamical degrees
\cite{Kho79,Tei79,Gro90}, this shows that
$h_\mu(F)>\log\Lambda$ forces all Lyapunov
exponents to be strictly positive.
For polynomial-like maps, the corresponding conclusions follow
from \cite[Theorem 4.1]{BR24}.
The local inverse branch construction
\cite{BDM08,BD19,BH24} 
and the Brin--Katok entropy formula
then give the analogue of Lemma \ref{lem:good-set}.
For the Bowen entropy statement, uniform expansion provides
the required branches.

\smallskip

The remaining potential-theoretic arguments also apply.
On $X$, compactness of normalized $\omega$-plurisubharmonic
functions \cite[Proposition 2.7]{GZ05} and smooth decreasing
approximation \cite[Theorem 1]{BKo07} give the analogue of
Lemma \ref{lem:holder-limit}.
By \cite[Theorem 7.2]{GZ05}, every locally pluripolar
set in $X$ is contained in the pole set of a global
$\omega$-psh function, so the integrability arguments
give the pluripolar conclusions.
The regularity consequences follow as in
Section \ref{subsec:restrictions}, using \cite{DN14}.

\bibliographystyle{amsalpha}
\bibliography{pluripolar_1111_V5}

\providecommand{\bysame}{\leavevmode\hbox to3em{\hrulefill}\thinspace}
\providecommand{\MR}{\relax\ifhmode\unskip\space\fi MR }
\providecommand{\MRhref}[2]{%
  \href{http://www.ams.org/mathscinet-getitem?mr=#1}{#2}
}
\providecommand{\href}[2]{#2}
\begin{thebibliography}{RHRHW26}

\bibitem[Abr59]{Abr59}
L.~M. Abramov, \emph{The entropy of a derived automorphism}, Dokl. Akad. Nauk
  SSSR \textbf{128} (1959), 647--650, In Russian.

\bibitem[AG01]{AG01}
David~H. Armitage and Stephen~J. Gardiner, \emph{Classical potential theory},
  Springer Monographs in Mathematics, Springer-Verlag, London, 2001.

\bibitem[Ale81]{Ale81}
Herbert Alexander, \emph{Projective capacity}, Recent developments in several
  complex variables, Ann. of Math. Stud., vol. 100, Princeton Univ. Press,
  Princeton, NJ, 1981, pp.~3--27.

\bibitem[BD01]{BriDuv01}
Jean-Yves Briend and Julien Duval, \emph{Deux caract{\'e}risations de la mesure
  d'{\'e}quilibre d'un endomorphisme de {$\mathbb{P}^k(\mathbb{C})$}}, Publ.
  Math. Inst. Hautes {\'E}tudes Sci. \textbf{93} (2001), 145--159.

\bibitem[BD19]{BD19}
Fran{\c{c}}ois Berteloot and Christophe Dupont, \emph{A distortion theorem for
  the iterated inverse branches of a holomorphic endomorphism of {$\mathbb
  P^k$}}, J. Lond. Math. Soc. (2) \textbf{99} (2019), no.~1, 153--172.

\bibitem[BD23]{BD23}
Fabrizio Bianchi and Tien-Cuong Dinh, \emph{Equilibrium states of endomorphisms
  of {$\mathbb P^k$} {I}: Existence and properties}, J. Math. Pures Appl. (9)
  \textbf{172} (2023), 164--201.

\bibitem[BD24]{BD24}
\bysame, \emph{Equilibrium states of endomorphisms of {$\mathbb{P}^k$}:
  Spectral stability and limit theorems}, Geom. Funct. Anal. \textbf{34}
  (2024), 1006--1051.

\bibitem[BDM08]{BDM08}
Fran{\c{c}}ois Berteloot, Christophe Dupont, and Laura Molino,
  \emph{Normalization of bundle holomorphic contractions and applications to
  dynamics}, Ann. Inst. Fourier (Grenoble) \textbf{58} (2008), no.~6,
  2137--2168.

\bibitem[BH24]{BH24}
Fabrizio Bianchi and Yan~Mary He, \emph{A {Ma{\~n}{\'e}--Manning} formula for
  expanding measures for endomorphisms of {$\mathbb P^k$}}, Trans. Amer. Math.
  Soc. \textbf{377} (2024), no.~11, 8179--8219.

\bibitem[BH26]{BH26Joining}
\bysame, \emph{Joining rigidity for rational maps}, 2026, arXiv:2609.14890.

\bibitem[BK83]{BK83}
Michael Brin and Anatole Katok, \emph{On local entropy}, Geometric Dynamics
  ({Rio de Janeiro}, 1981), Lecture Notes in Math., vol. 1007, Springer,
  Berlin, 1983, pp.~30--38.

\bibitem[BK07]{BKo07}
Zbigniew B{\l}ocki and S{\l}awomir Ko{\l}odziej, \emph{On regularization of
  plurisubharmonic functions on manifolds}, Proc. Amer. Math. Soc. \textbf{135}
  (2007), no.~7, 2089--2093.

\bibitem[BLS93]{BLS93}
Eric Bedford, Mikhail Lyubich, and John Smillie, \emph{Polynomial
  diffeomorphisms of {$\mathbb C^2$}. {IV}: The measure of maximal entropy and
  laminar currents}, Invent. Math. \textbf{112} (1993), 77--125.

\bibitem[BR24]{BR24}
Fabrizio Bianchi and Karim Rakhimov, \emph{Strong probabilistic stability in
  holomorphic families of endomorphisms of {$\mathbb{P}^k(\mathbb{C})$} and
  polynomial-like maps}, Int. Math. Res. Not. IMRN \textbf{2024} (2024),
  no.~12, 9875--9892.

\bibitem[BS91]{BS91}
Eric Bedford and John Smillie, \emph{Polynomial diffeomorphisms of {$\mathbb
  C^2$}: currents, equilibrium measure and hyperbolicity}, Invent. Math.
  \textbf{103} (1991), 69--99.

\bibitem[BT82]{BT82}
Eric Bedford and B.~A. Taylor, \emph{A new capacity for plurisubharmonic
  functions}, Acta Math. \textbf{149} (1982), 1--40.

\bibitem[Can14]{Can14}
Serge Cantat, \emph{Dynamics of automorphisms of compact complex surfaces},
  Frontiers in Complex Dynamics, Princeton Math. Ser., vol.~51, Princeton
  University Press, Princeton, NJ, 2014, pp.~463--514.

\bibitem[CFS82]{CFS82}
I.~P. Cornfeld, S.~V. Fomin, and Ya.~G. Sinai, \emph{Ergodic theory},
  Grundlehren Math. Wiss., vol. 245, Springer-Verlag, New York, 1982.

\bibitem[DN14]{DN14}
Tien-Cuong Dinh and Vi{\^e}t-Anh Nguy{\^e}n, \emph{Characterization of
  {Monge--Amp\`ere} measures with {H{\"o}lder} continuous potentials}, J.
  Funct. Anal. \textbf{266} (2014), no.~1, 67--84.

\bibitem[DNS10]{DNS10}
Tien-Cuong Dinh, Vi{\^e}t-Anh Nguy{\^e}n, and Nessim Sibony, \emph{Exponential
  estimates for plurisubharmonic functions and stochastic dynamics}, J.
  Differential Geom. \textbf{84} (2010), no.~3, 465--488.

\bibitem[DR20]{DR20}
Christophe Dupont and Axel Rogue, \emph{Dimension of ergodic measures and
  currents on {$\mathbb{CP}(2)$}}, Ergodic Theory Dynam. Systems \textbf{40}
  (2020), no.~8, 2131--2155.

\bibitem[DS03]{DS03}
Tien-Cuong Dinh and Nessim Sibony, \emph{Dynamique des applications d'allure
  polynomiale}, J. Math. Pures Appl. (9) \textbf{82} (2003), no.~4, 367--423.

\bibitem[DS09]{DS09SP}
\bysame, \emph{Super-potentials of positive closed currents, intersection
  theory and dynamics}, Acta Math. \textbf{203} (2009), no.~1, 1--82.

\bibitem[DS10a]{DS10}
\bysame, \emph{Dynamics in several complex variables: endomorphisms of
  projective spaces and polynomial-like mappings}, Holomorphic Dynamical
  Systems, Lecture Notes in Math., vol. 1998, Springer, Berlin, 2010,
  pp.~165--294.

\bibitem[DS10b]{DS10Aut}
\bysame, \emph{Super-potentials for currents on compact {K\"ahler} manifolds
  and dynamics of automorphisms}, J. Algebraic Geom. \textbf{19} (2010), no.~3,
  473--529.

\bibitem[dT08]{dT08}
Henry de~Th{\'e}lin, \emph{Sur les exposants de {Lyapounov} des applications
  m{\'e}romorphes}, Invent. Math. \textbf{172} (2008), no.~1, 89--116.

\bibitem[dTV15]{dTV15}
Henry de~Th{\'e}lin and Gabriel Vigny, \emph{On the measures of large entropy
  on a positive closed current}, Math. Z. \textbf{280} (2015), no.~3--4,
  919--944.

\bibitem[Duj04]{Duj04HL}
Romain Dujardin, \emph{{H{\'e}non}-like mappings in {$\mathbb C^2$}}, Amer. J.
  Math. \textbf{126} (2004), no.~2, 439--472.

\bibitem[Duj05]{Duj05}
\bysame, \emph{Structure properties of laminar currents on {$\mathbb P^2$}}, J.
  Geom. Anal. \textbf{15} (2005), no.~1, 25--47.

\bibitem[Dup12]{Dup12}
Christophe Dupont, \emph{Large entropy measures for endomorphisms of
  {$\mathbb{CP}^k$}}, Israel J. Math. \textbf{192} (2012), no.~2, 505--533.

\bibitem[Fel71]{Fel71}
William Feller, \emph{An introduction to probability theory and its
  applications}, 2 ed., vol.~II, John Wiley \& Sons, New York, 1971.

\bibitem[FH12]{FH12}
De-Jun Feng and Wen Huang, \emph{Variational principles for topological
  entropies of subsets}, J. Funct. Anal. \textbf{263} (2012), no.~8,
  2228--2254.

\bibitem[FLM83]{FLM83}
Alexandre Freire, Artur Lopes, and Ricardo Ma{\~n}{\'e}, \emph{An invariant
  measure for rational maps}, Bol. Soc. Brasil. Mat. \textbf{14} (1983), no.~1,
  45--62.

\bibitem[FS95]{FS95}
John~Erik Forn{\ae}ss and Nessim Sibony, \emph{Complex dynamics in higher
  dimension. {II}}, Modern Methods in Complex Analysis ({Princeton}, {NJ},
  1992), Ann. of Math. Stud., vol. 137, Princeton University Press, Princeton,
  NJ, 1995, pp.~135--182.

\bibitem[Gro90]{Gro90}
Mikha{\"i}l Gromov, \emph{Convex sets and {K{\"a}hler} manifolds}, Advances in
  Differential Geometry and Topology (Franco Tricerri, ed.), World Scientific,
  Teaneck, NJ, 1990, pp.~1--38.

\bibitem[Gro03]{Gro03}
\bysame, \emph{On the entropy of holomorphic maps}, Enseign. Math. (2)
  \textbf{49} (2003), no.~3--4, 217--235.

\bibitem[GZ05]{GZ05}
Vincent Guedj and Ahmed Zeriahi, \emph{Intrinsic capacities on compact
  {K{\"a}hler} manifolds}, J. Geom. Anal. \textbf{15} (2005), no.~4, 607--639.

\bibitem[Kho79]{Kho79}
Askold~G. Khovanskii, \emph{The geometry of convex polyhedra and algebraic
  geometry}, Uspekhi Mat. Nauk \textbf{34} (1979), no.~4, 160--161, In Russian.

\bibitem[LY85]{LY85}
Fran{\c{c}}ois Ledrappier and Lai-Sang Young, \emph{The metric entropy of
  diffeomorphisms. {Part II}: Relations between entropy, exponents and
  dimension}, Ann. of Math. (2) \textbf{122} (1985), no.~3, 540--574.

\bibitem[Lyu83]{Lyu83}
Mikhail Lyubich, \emph{Entropy properties of rational endomorphisms of the
  {Riemann} sphere}, Ergodic Theory Dynam. Systems \textbf{3} (1983), no.~3,
  351--385.

\bibitem[Ma{\~n}83]{Man83}
Ricardo Ma{\~n}{\'e}, \emph{On the uniqueness of the maximizing measure for
  rational maps}, Bol. Soc. Brasil. Mat. \textbf{14} (1983), no.~1, 27--43.

\bibitem[Ma{\~n}88]{Man88}
\bysame, \emph{The {Hausdorff} dimension of invariant probabilities of rational
  maps}, Dynamical Systems ({Valparaiso}, 1986), Lecture Notes in Math., vol.
  1331, Springer, Berlin, 1988, pp.~86--117.

\bibitem[Nor97]{Nor97}
J.~R. Norris, \emph{Markov chains}, Cambridge Series in Statistical and
  Probabilistic Mathematics, vol.~2, Cambridge University Press, Cambridge,
  1997.

\bibitem[RHRHW26]{RHW26}
Federico Rodriguez~Hertz, Jana Rodriguez~Hertz, and Weisheng Wu,
  \emph{Ergodicity of the geodesic flow on surfaces of nonpositive curvature},
  2026, arXiv:2610.27697.

\bibitem[Sko72]{Sko72}
Henri Skoda, \emph{Sous-ensembles analytiques d'ordre fini ou infini dans
  {$\mathbb{C}^n$}}, Bull. Soc. Math. France \textbf{100} (1972), 353--408.

\bibitem[Tei79]{Tei79}
Bernard Teissier, \emph{Du th{\'e}or{\`e}me de l'index de {Hodge} aux
  in{\'e}galit{\'e}s isop{\'e}rim{\'e}triques}, C. R. Acad. Sci. Paris S{\'e}r.
  A-B \textbf{288} (1979), no.~4, A287--A289.

\bibitem[UZ13]{UZ13}
Mariusz Urba{\'n}ski and Anna Zdunik, \emph{Equilibrium measures for
  holomorphic endomorphisms of complex projective spaces}, Fund. Math.
  \textbf{220} (2013), no.~1, 23--69.

\end{thebibliography}

\end{document}